\documentclass[sn-mathphys]{sn-jnl}
\usepackage[numbers,square]{natbib}
\usepackage{graphicx}%
\usepackage{multirow}%
\usepackage{amsmath,amssymb,amsfonts}%
\usepackage{amsthm}%
\usepackage{mathrsfs}%
\usepackage[title]{appendix}%
\usepackage{xcolor}%
\usepackage{textcomp}%
\usepackage{manyfoot}%
\usepackage{booktabs}%
\usepackage{algorithm}%
\usepackage{algorithmicx}%
\usepackage{algpseudocode}%
\usepackage{listings}%
\usepackage{mathtools}

\usepackage{tikz}
\usetikzlibrary{arrows}
\usepackage{subcaption}

\makeatletter
\def\moverlay{\mathpalette\mov@rlay}
\def\mov@rlay#1#2{\leavevmode\vtop{%
   \baselineskip\z@skip \lineskiplimit-\maxdimen
   \ialign{\hfil$\m@th#1##$\hfil\cr#2\crcr}}}
\newcommand{\charfusion}[3][\mathord]{
    #1{\ifx#1\mathop\vphantom{#2}\fi
        \mathpalette\mov@rlay{#2\cr#3}
      }
    \ifx#1\mathop\expandafter\displaylimits\fi}
\makeatother

\theoremstyle{thmstyleone}%
\newtheorem{theorem}{Theorem}%  meant for continuous numbers
\newtheorem{lemma}[theorem]{Lemma}%
\newtheorem{observation}[theorem]{Observation}%
\newtheorem{corollary}[theorem]{Corollary}%

\theoremstyle{thmstyletwo}%
\newtheorem{example}{Example}%
\newtheorem{remark}{Remark}%

\theoremstyle{thmstylethree}%
\newtheorem{definition}{Definition}%

\begin{document}

%% ================================

\title[Directed Transit Functions]{\textbf{Directed Transit Functions}}

\newcommand{\AX}[1]{\textsf{\upshape{(#1)}}}

\newcommand{\TODO}[1]{\begingroup\color{red}#1\endgroup}
\newcommand{\PFS}[1]{\begingroup\color{blue}#1\endgroup}

\author[1]{\fnm{Arun} \sur{Anil}}\email{arunanil93@gmail.com}
\author[1]{\fnm{Manoj}
  \sur{Changat}}\email{mchangat@keralauniversity.ac.in}
\author[1]{\fnm{Lekshmi Kamal}
  \sur{K-Sheela}}\email{lekshmisanthoshgr@gmail.com}
\author[1]{\fnm{Ameera} \sur{Vaheeda Shanavas}}\email{ameerasv@gmail.com}
\author[2,3]{\fnm{John J.} \sur{Chavara}}\email{jjchavara@gmail.com}
\author[4,5]{\fnm{Prasanth G.}
  \sur{Narasimha-Shenoi}}\email{prasanthgns@gmail.com}
\author[6,7]{\fnm{Bruno J.}
  \sur{Schmidt}}\email{bruno@bioinf.uni-leipzig.de}
\author[7,8,6,9,10,11]{\fnm{Peter F.}
  \sur{Stadler}}\email{studla@bioinf.uni-leipzig.de}

\affil[1]{\orgdiv{Department of Futures Studies}, \orgname{University of
    Kerala}, \orgaddress{\street{Karyavattom Campus}
    \postcode{695 581} \city{Thiruvananthapuram}, 
    \country{India}}}
 
\affil[2]{\orgdiv{Department of Mathematics}, \orgname{University of
    Kerala},  \orgaddress{\street{Karyavattom Campus}
    \postcode{695 581} \city{Thiruvananthapuram}, 
    \country{India}}}

\affil[3]{\orgdiv{Department of Mathematics}, \orgname{St.\ Berchmans College},
  \orgaddress{\street{Changanassery}, \city{Kottayam}, \postcode{686101},
    \country{India}}}

\affil[4]{\orgdiv{Department of Mathematics}, \orgname{Government College},
  \orgaddress{\street{Chittur}, \city{Palakkad}, \postcode{678 104},
    \country{India}}}

\affil[5]{\orgdiv{Department of Collegiate Education},
  \orgname{Government of Kerala}, \orgaddress{\city{Thiruvananthapuram},
    \country{India}}}

\affil[6]{\orgname{Max Planck Institute for Mathematics in the Sciences},
  \orgaddress{\street{Inselstra{\ss}e 22},
    \postcode{D-04103} \city{Leipzig}, \country{Germany}}}

\affil[7]{\orgdiv{Bioinformatics Group, Department of Computer Science \&
    Interdisciplinary Center for Bioinformatics, \orgname{Leipzig
      University}, \orgaddress{\street{H{\"a}rtelstra{\ss}e 16–18},
      \postcode{D-04107} \city{Leipzig}, \country{Germany}}}

\affil[8]{\orgdiv{School for Embedded and Composite Artificial
    Intelligence (SECAI)},
  \orgname{Leipzig
    University}, \city{Leipzig}, \country{Germany}}}

\affil[9]{\orgdiv{Department of Theoretical Chemistry}, \orgname{University
    of Vienna}, \orgaddress{\street{W{\"a}hringerstra{\ss}e 17},
    \postcode{A-1090} \city{Wien}, \country{Austria}}}

\affil[10]{\orgdiv{Facultad de Ciencias}, \orgname{Universidad National de
    Colombia}, \orgaddress{\city{Bogot{\'a}}, \country{Colombia}}}

\affil[11]{\orgname{Santa Fe Institute}, \orgaddress{\street{1399 Hyde Park
      Rd.}, \city{Santa Fe}, \state{NM} \postcode{87501}, \country{USA}}}

\abstract{Transit functions were introduced as models of betweenness on
  undirected structures. Here we introduce directed transit function as the
  directed analogue on directed structures such as posets and directed
  graphs. We first show that betweenness in posets can be expressed by
  means of a simple set of first order axioms. Similar characterizations
  can be obtained for graphs with natural partial orders, in particular,
  forests, trees, and mangroves. Relaxing the acyclicity conditions leads
  to a generalization of the well-known geometric transit function to the
  directed structures. Moreover, we discuss some properties of the directed
  analogues of prominent transit functions, including the all-paths,
  induced paths, and shortest paths (or interval) transit functions.
  Finally we point out some open questions and directions for future work.
}

\title{Directed Transit Functions}


%\classnbr{05C38, 05C69, 05C99} 

\maketitle

\section{Introduction}

Transit functions $R:V\times V\to 2^V$ on a non-empty finite set $V$ were
introduced as an abstract model of betweenness, see \cite{Mulder:08a}. In
the most general form, $R$ satisfies the three basic axioms (t1)
$\{u,v\}\in R(u,v)$, (t2) $R(u,v)=R(v,u)$, and (t3)
$R(u,u)=\{u\}$. Important examples of transit functions are obtained by
defining $R(u,v)$ as the set of all vertices that lie on a certain class of
paths and walks connecting $u$ and $v$. This includes shortest paths
\cite{mu-80,nebe-94,mune-09}, induced
paths\cite{nebe-02,Changat:04i,mcjmhm-10}, all-paths \cite{Changat:01},
toll walks \cite{toll1,lcp}, and weak toll walks \cite{weaktoll} in
graphs. In a very different context, transit functions also appear as
descriptions of so-called clustering systems, in which each cluster $C$ is
the intersection of all clusters that contain a pair of suitably chosen
points $x,y\in C$
\cite{Bandelt:89,Barthelemy:08,Bertrand:24,Changat:18a,Changat:23x,Shanavas:24a}.

Recently, axioms (t1) and (t3) were relaxed in \cite{Changat:17} to study
the interval function of not necessarily connected undirected graphs, using the axioms (t1*) $R(u,v)\ne\emptyset$ implies $u\in R(u,v)$, and
(compt) $R(u,v)\ne\emptyset$ and $R(v,w)\ne\emptyset$ then
$R(u,w)\ne\emptyset$.  This generalization does not seem to have drawn much
attention so far beyond the work reported in \cite{Changat:17}.

Interval functions satisfying (t1), (t2), and (t3) have been constructed
for directed graphs (digraphs) in \cite{Chartrand:00} as the union of all
shortest directed path from $u$ to $v$ and $v$ to $u$, setting
$R(u,v)=\{u,v\}$ if no such path exists. This notion has been explored in
particular for oriented (antisymmetric) graphs in relation to convexities
derived from these intervals. Here, however, we will pursue a different
avenue, considering also the intervals as directed structures.

Measures of betweenness, in particular \emph{betweenness centrality}, i.e.,
the average fraction of shortest $s$-$t$-paths passing through a vertex
$v$, have been studied extensively also for directed networks
\cite{White:94,Brandes:01}.  Surprisingly, the theory of transit functions,
so far has not been extended to directed structures.

Notions of paths connecting two points are not restricted to graphs.
Directed hypergraphs have received increasing attention of constructive
dynamical systems including chemical reaction networks
\cite{Estrada:11,Fagerberg:18a}. As in graphs, hyperpaths are defined as
alternating sequences of vertices and edges
$(x_1,E_1,x_2,\dots,x_{n-1},E_n,x_n)$ subject to a wide variety of
additional conditions, see
e.g.\ \cite{Thakur:09,Ausiello:12,Dharmarajan:15}.

A directed form of betweenness also arises from bubble structures in
directed graphs, defined as sets of vertices reachable from an entrance $s$
and from which an exit $t$ can be reached. Such structures appear in
conjunction with additional isolation properties e.g.\ in the context of
genome assembly \cite{Onodera:13,Paten:18,Gaertner:18b} and as a model of
alternative splicing \cite{Sammeth:09} in computational biology.

In this contribution we introduce transit functions for directed
structures. As a very natural starting point we will first consider partial
orders in Sect.~\ref{sect:PO}. While the standard notion of betwenness for
partial orders using symmetric transit functions, cannot be characterized
by a single first-order sentence \cite{Courcelle:20}, its directed
counterpart turns out to have a simple and very natural
characterization. We explore the directed transit functions on posets
further by describing general directed acyclic graphs (DAG) as well as some special subclasses such
as rooted trees in more detail. In Sect.~\ref{sect:reach} we relax our
assumptions and consider directed betweenness in general reachability
relations. This leads us in particular to a natural generalization of
Nebesk{\'y}'s \emph{geometric transit functions} \cite{Nebesky:01}. We then
briefly discuss constructions of directed transit functions on graphs that
have been studied in detail in the undirected case, namely the all-paths,
the induced-path, and shortest-path transit functions. As we shall see,
there are non-trivial differences between undirected and directed transit
functions.

\section{Notation}

We consider finite, loopless, simple directed graphs (digraphs) $G=(V,E)$
and write $V(G)$ and $E(G)$ to refer explicitly to the vertex and edge set
of $G$. A directed edge from $u$ to $v$ in $G$ will be denoted by $(u,v)\in
E(G)$. A walk in $G$ is an alternating sequence $x_0,e_1,x_1,\dots,e_n,x_n$
of vertices $x_i\in V(G)$ and edges $e_i=(x_{i-1},x_i)\in E(G)$. Since the
edges are implied by the vertices, we will simply write
$(x_0,x_1,\dots,x_n)$ for the walk starting from vertex $x_0$ and ending in
the vertex $x_n$. The number $n$ of edges is the length of a walk. A path
is a walk such that $x_i\ne x_j$ for all $i\ne j$. A directed cycle is
composed of a path $(x_1,\dots,x_n)$ and the edge $(x_n,x_1)$. 

A digraph is symmetric if $(u,v)\in E(G)$ implies $(v,u)\in
E(G)$. Symmetric digraphs can also be interpreted as undirected graphs.  A
digraph $G$ is {strongly connected} if for any two vertices $u,v\in V(G)$
there exists a path from $u$ to $v$ and from $v$ to $u$. It is unilaterally
connected if for any two vertices $u,v\in V(G)$ there exists a path from
$u$ to $v$ or from $v$ to $u$. Finally, $G$ is weakly connected if for any
two vertices $u,v\in V(G)$ there is a sequence $(u=x_0,x_1,\dots,x_n=v)$
such that for $1\le i\le n$ we have $(x_{i-1},x_i)\in E(G)$ or
$(x_{i},x_{i-1})\in E(G)$.

A vertex $x\in V(G)$ is a source if it has no incoming edges
$(u,x)\in E(G)$, and a sink if it has no outgoing edges $(x,v)\in E(G)$.
If $G$ is a DAG, a vertex $x\in V(G)$ is called a hybrid vertex if it has
more than one incoming edge. A source vertex $w\in V(G)$ is called root
vertex, or simply root, if there is a path from $w$ to every other vertex
in $V(G)$. A DAG is rooted (at $w$) if it contains a root vertex $w$. A
graph containing a root vertex is always weakly connected by definition.
If $(x,y)\in E(G)$ in a DAG $G$, then $y$ is called a successor of $x$, and
$x$ a predecessor of $y$.

$H=(W,F)$ is a sub(di)graph of $G=(V,E)$ if $H$ is a digraph,
$W\subseteq V$ and $F\subseteq E$. The subgraph $H=G[W]$ is induced by $W$
if for all $u,v\in W$ we have $(u,v)\in E$ if and only if $(u,v)\in F$. A
strongly (or weakly) connected component of $G$ is a maximal subset
$W\subseteq V(G)$ such that the subgraph $G[W]$ induced by $W$ is strongly
(or weakly) connected. A rooted (directed) tree is a rooted DAG without
hybrid vertices. Thus any path between two vertices is unique. A (directed)
forest $\vec{F}$ is the disjoint union of rooted trees.

A partial order is a reflexive, antisymmetric, transitive relation on a
non-empty set $V$, i.e., (i) $x\succeq x$ for all $x\in V$, (ii)
$x\succeq y$ and $y\succeq x$ implies $x=y$ for all $x,y\in V$, and (iii)
$x\succeq z$ and $z\succeq y$ implies $x\succeq y$ for all $x,y,z\in V$.

\section{Directed Transit Functions on Partially Ordered Sets}
\label{sect:PO}

\subsection{Motivation}

Partially ordered sets (posets) are probably the simplest ``directed''
structures, comprising a set $V$ endowed with a reflexive, transitive, and
antisymmetric relation. Posets are associated with a straightforward notion
of betweenness, which can be captured as follows:
\begin{definition}
  Let $(V,\succeq)$ be a non-empty finite partial order. Then the function
  $R_{\succeq}: V\times V\to 2^V$  define by 
  \begin{equation}
  \label{eq:poset}
  R_{\succeq}(u,v) \coloneqq \{x\in V| u\succeq x\succeq v\}
  \end{equation}
  is the \emph{directed transit function} of $\succeq$.
\end{definition}
For incomparable elements $u$ and $v$, i.e., for $u$ and $v$ such that
$u\not\succeq v$ and $v\not\succeq u$, we use $u \parallel v$.  Note that
$R_{\succeq}(u,v)=\emptyset$ if and only if $u\not\succeq v$.  Moreover, we
have $u\parallel v$ if and only if
$R_{\succeq}(u,v)=R_{\succeq}(v,u)=\emptyset$.
  
We will usually simply write $R$ instead of $R_{\succeq}$ when there is no
danger of confusion. We note that \emph{symmetric} transit functions also
have been associated with posets \cite{Mulder:08a}. Similar to
eqn.(\ref{eq:poset}), if $u$ and $v$ are comparable, $T(u,v)=\{x\in
V| u\succeq x\succeq v\}$ if $u\succeq v$ and $T(v,u)=\{x\in V| v\succeq
x\succeq u\}$ if $v\succeq u$. In contrast to our definition, however,
$T(u,v)=\{u,v\}$ is used if $u$ and $v$ are incomparable. We use the empty
set here in the same manner as in the generalized transit functions on
disconnected graphs in \cite{Changat:17}.

\begin{lemma}
  \label{lem:poset1}
  Let $R$ be the transit function of a partial order. Then $R$ satisfies
  the following statements:
  \begin{itemize}
  \item[\AX{t1}]   If $R(u,v)\ne\emptyset$ then $\{u,v\}\subseteq R(u,v)$
    for all $u,v\in V$.
  \item[\AX{t2a}]  If $R(u,v)\ne\emptyset$ and $u\ne v$ then
    $R(v,u)=\emptyset$ for all $u,v\in V$.
  \item[\AX{t3}]   $R(u,u)=\{u\}$ for all $u\in V$. 
  \item[\AX{tr1}] If $R(u,w)\ne \emptyset$ and $R(w,v)\ne \emptyset$ then
    $R(u,w)\cup R(w,v)\subseteq R(u,v)$ for all $u,v\in V$ with $u\ne v$.
  \item[\AX{tr2}] If $R(u,w)=\emptyset$ or $R(w,v)=\emptyset$ then $w\notin
    R(u,v)$ for all $u,v\in V$ with $u\ne v$.
  \end{itemize}
\end{lemma}
\begin{proof}
  \AX{t1} follows directly from the definition, and $u\succeq u$ yields
  $R(u,u)=\{u\}$ for all $u\in V$, i.e., \AX{t3}.  Since $R(u,v)\ne\emptyset$
  implies $u\succeq v$ and $u\ne v$, we have $v\not\succeq u$ and thus
  $R(u,v)=\emptyset$, and consequently \AX{t2a} is satisfied. Similarly,
  $R(u,v)\ne \emptyset$ and $R(v,w)\ne \emptyset$ implies $u\succeq v$ and
  $v\succeq w$, which by transitivity of the partial order yields $u\succeq
  x\succeq v\succeq w$ and thus equ.(\ref{eq:poset}) implies
  $R(u,v)\subseteq R(u,w)$. Analogously, we obtain $R(v,w)\subseteq
  R(u,w)$, and hence \AX{tr1} is satisfied. Finally, $u\not\succeq x$ or
  $x\not\succeq v$ implies immediately that $x\notin R(u,v)$, and thus
  \AX{tr2} holds.
\end{proof}
We note that in the presence of \AX{t2a}, the statements of \AX{tr1} and
\AX{tr2} trivially extends to the excluded cases $u=v$ as an immediate
consequence of \AX{t3}. We explicitly exclude $u=v$ in the statement of the
axioms \AX{tr1} and \AX{tr2} because this will be relevant later in the
context of reachability relations.

Before we proceed we consider two simple consequences of \AX{tr1}.
\begin{observation}
  \label{obs:tr0}
  Let $R: V\times V\to 2^V$ be a function satisfying \AX{t3}, and
  \AX{tr1}. Then $R$ satisfies
  \begin{itemize}
  \item[\AX{tr0}] If $u\ne v$, $R(u,w)\ne \emptyset$ and $R(w,v)\ne
    \emptyset$ then $w\in R(u,v)$.
  \item[\AX{t0}] If $R(u,w)\ne \emptyset$ and $R(w,v)\ne \emptyset$ then
    $R(u,v)\ne\emptyset$. 
  \end{itemize}
\end{observation}
\begin{proof}
  If $u\ne v$, and $R(u,w)\ne \emptyset$ and $R(w,v)\ne \emptyset$ then
  $w\in R(u,w)\cup R(w,v)\subseteq R(u,v)\ne\emptyset$, and thus \AX{tr0}
  holds and the implication of \AX{t0} is true.  If $u=v$, then \AX{t3}
  yields $u\in R(u,v)\ne\emptyset$, and hence \AX{t0} holds for all $u,v\in
  V$.
\end{proof}

We note that under the assumptions of Observation~\ref{obs:tr0}, \AX{tr0}
implies \AX{t0}, while the converse is not true.  In \cite{Changat:17},
\AX{t0} (there denoted by ``compt''), \AX{t1}, and \AX{t3} appear in this
form for undirected, but not necessarily connected, graphs. The analogy of
the basic properties of $R$ in Lemma~\ref{lem:poset1} on properties with
symmetric transit functions is to be motivated by the following
\begin{definition}
  Let $V$ be a finite set. A \emph{directed transit function} is a function
  $R:V\times V\to 2^V$ satisfying \AX{t0}, \AX{t1} and \AX{t3}.
\end{definition} 

\subsection{Characterization of the Directed Transit Functions of Posets} 

Given a directed transit function, we define a relation $\succeq_R$ by 
$u\succeq_R v$ if and only if $R(u,v)\ne\emptyset$. Next, we show that the 
construction of $R_{\succeq}$ preserves the partial order in the following 
sense:
\begin{lemma}
  Let $\succeq$ be a partial order and $R_{\succeq}$ its directed transit
  function. Then $\succeq_{R_{\succeq}}=\succeq$.
\end{lemma}
\begin{proof}
  If $u\succeq v$ then $R_{\succeq}(u,v)\ne\emptyset$ and thus
  $u\succeq_{R_{\succeq}} v$. If $u\not\succeq v$ then
  $R_{\succeq}(u,v)=\emptyset$ and thus  $u\not\succeq_{R_{\succeq}} v$.
  Hence $\succeq_{R_{\succeq}}=\succeq$.
\end{proof}

\begin{definition} 
  A directed transit function is a \emph{poset function} if it is
  antisymmetric, \AX{t2a}, and satisfies \AX{tr1} and \AX{tr2}.
\end{definition}

The main result of this section is that poset functions faithfully describe
posets.
\begin{theorem}
  Let $R$ be a poset function.  Then the relation $\succeq_R$ defined by
  $u\succeq_R v$ if and only if $R(u,v)\ne\emptyset$ is a partial order.
  Moreover $R_{\succeq_R}=R$.
\end{theorem}
\begin{proof}
  First we note that by \AX{t3} we have $u\succeq_R u$ for all $u\in V$.
  Hence $\succeq_R$ is reflexive. Suppose $u\ne v$. Then $u\succeq_R v$
  implies $R(u,v)\ne\emptyset$ and thus $R(v,u)=\emptyset$ by \AX{t2a},
  i.e., $v\not\succeq u$, i.e., $\succeq_R$ is antisymmetric. Finally
  assume that $u\succeq_R v$ and $v\succeq_R w$. Thus
  $R(u,v)\ne\emptyset$ and $R(v,w)\ne\emptyset$, whence \AX{tr1} implies
  $R(u,w)\ne\emptyset$, and thus $u\succeq_R w$, i.e, $\succeq_R$ is
  transitive. Thus $\succeq_R$ is a partial order.

  By definition, $x\in R_{\succeq_R}(u,v)$ if and only if $u\succeq_R x$ and $x
  \succeq_R v$, which is the case if and only if $R(u,x)\ne\emptyset$ and
  $R(x,v)\ne\emptyset$; thus \AX{tr1}, which implies $x\in R(u,v)$.  If
  $x\notin R_{\succeq_R}(u,v)$ then $u\not\succeq_R x$ or $x\not\succeq_R
  v$, i.e., $R(u,x)=\emptyset$ or $R(x,v)=\emptyset$, which by \AX{tr2}
  implies $x\notin R(u,v)$. Hence $R{\succeq_R}=R$.
\end{proof}

Before we proceed we note that \AX{tr1} and \AX{tr2} can be expressed in
the following, more compact form.
\begin{lemma}
  Let $R$ be a directed transit function. Then $R$ satisfies \AX{tr1} and
  \AX{tr2} if and only $R$ satisfies, for all $u,v,w\in V$ with $u\ne v$:
  \begin{itemize}
    \item[\AX{tr}] $R(u,w)\ne\emptyset$ and $R(w,v)\ne\emptyset$ if and
      only if $w\in R(u,v)$
  \end{itemize}
\end{lemma}
We remark that in the presence of \AX{t3} and \AX{t2a}, the equivalence
also remains valid for $u=v$.
\begin{proof}
  Assume \AX{tr1}, \AX{tr2}. If $w\in R(u,v)$ and $R(u,w)=\emptyset$ then
  $w\notin R(u,v)$ by \AX{tr2}, a contradiction. The same argument can be
  repeated for $R(w,v)=\emptyset$, hence the if-part of \AX{tr} hold.  For
  the only-if part, assume $R(u,w)\neq\emptyset$ and $R(w,v)\neq\emptyset$
  and $w\notin R(u,v)$.  By \AX{tr1}, we have $R(u,w)\cup R(w,v)\subseteq
  R(u,v)$ and since $w\in R(u,w),R(w,v)$ by \AX{t1} we have $w\in R(u,v)$,
  a contradiction. Thus \AX{tr1} and \AX{tr2} imply \AX{tr}.
  
  For the converse, assume \AX{tr}. For \AX{tr1}, assume
  $R(u,w)\neq\emptyset$ and $R(w,v)\neq\emptyset$ and $R(u,w)\cup
  R(w,v)\not\subseteq R(u,v)$. By \AX{tr}, we know $w\in R(u,v)$ and by
  \AX{t1}, that $u,v\in R(u,v)$, hence there must be an $x\in R(u,w)$ or
  $x\in R(w,v)$ such that $x\notin R(u,v)$. If $x\in R(u,w)$ then
  $R(u,x)\ne\emptyset$, $R(x,w)\ne\emptyset$, and $R(w,v)\neq\emptyset$;
  Now $x\in R(x,v)$ and \AX{tr} implies $x\in R(u,v)$, a contradiction.  We
  argue analogously for $x\in R(w,u)$.  \AX{tr2} is the contrapositive of
  the only-if direction in \AX{tr}.
\end{proof}
Since a poset function satisfies \AX{tr1} and \AX{tr2} we have
\begin{remark}\label{remark_poset_tr}
  A poset function $R$ always satisfies \AX{tr}
\end{remark}
The following property of poset functions will be convenient for the proofs
in the following sections.
\begin{lemma}\label{lemma-axiom-q}
  Let $R$ be a poset function. Then $R$ satisfies, for all $u,v,w\in V$ 
  \begin{itemize}
  \item[\AX{q}] $R(u,w)\cap R(w,v)\in \{w,\emptyset \}$
  \end{itemize}
\end{lemma}
\begin{proof} 
  \par\noindent \AX{q}.\quad
  Suppose there is $x\in R(u,w)\cap R(w,v)$ with
  $x\ne w$.  Since $x\in R(u,w)$ and $x\in R(w,v)$, the contrapositive of
  \AX{tr2} implies $R(x,w)\ne\emptyset$ and $R(w,x)\ne\emptyset$,
  contradicting \AX{t2a}. Thus, if $R(u,w)\ne\emptyset$ and
  $R(w,v)\ne\emptyset$ then $R(u,w)\cap R(w,v)=\{w\}$. Clearly, $R(u,w)\cap
  R(w,v)=\emptyset$ if (and only if) $R(u,w)=\emptyset$ and
  $R(w,v)=\emptyset$.
\end{proof}

\subsection{Directed Acyclic Graphs}

A graph $G=(V,E)$ is acyclic if it does not contain a directed cycle. The
directed paths in directed acyclic graph (DAG) define a partial order
$\succeq$ on the vertex set $V$ such that $u\succeq v$ if and only if there
is a directed path from $u$ to $v$. Thus the directed all-path transit
function $R_G$ of $G$ coincides with $R_{\succeq_G}$. This begs the
questions whether the DAG $G$ can be identified from $R_G$.

A key concept linking transit functions with graphs is the \emph{underlying
graph} \cite{Mulder:08a}. The idea directly carries over to the directed
case:
\begin{definition}
  Let $R$ be a directed transit function on $V$. The underderlying digraph
  of $R$ is $G_R=(X,E_R)$, where $(x,y)\in E_R$ if and only if
  $R(x,y)=\{x,y\}$ and $x\neq y$.
\end{definition}
The underlying graph $G_R$ is symmetric if $R$ is symmetric (i.e., if $R$
satisfies (t2) $R(u,v)=R(v,u)$.  In this case $G_R$ corresponds to the
undirected underlying graph defined for symmetric transit functions.

\begin{lemma}
  \label{lem:upathT1}
  Let $R$ be a poset function.  If $R(u,v)\ne\emptyset$ and $u\ne v$ then
  there are $x,y\in R(u,v)$ such that $R(u,x)=\{u,x\}$ and
  $R(y,v)=\{y,v\}$.
\end{lemma}
\begin{proof}
  If $R(u,v)=\{u,v\}$ there is nothing to show. Otherwise, there is $w\in
  R(u,v)\setminus\{u,v\}$, \AX{q} implies $R(u,w)\cap R(w,v)=\{w\}$, and
  thus $u\notin R(w,v)$ and $v\notin R(u,w)$. On the other hand, \AX{tr1}
  ensures $R(u,w)\subseteq R(u,v)$ and $R(w,v)\subseteq R(u,v)$, and thus
  $|R(u,w)|,|R(w,v)|< |R(u,v)|$. Applying the argument separately to
  $R(u,w)$ and $R(w,v)$ yields, elements $x_2\in R(u,w)$, $y_2\in R(w,v)$,
  with $R(u,x_2)\subsetneq R(u,w)$ and, $R(y_2,v)\subseteq R(w,v)$ and,
  thus a sequence of transit sets having elements $x_{k-1},x_{k},
  x_{k+1},\ldots$, and $y_{\ell-1},y_{\ell}, y_{\ell+1},\ldots$, with
  $x_k\in R(u,x_{k-1})$, $x_{k+1}\in R(u,x_{k})$ and $y_{\ell}\in
  R(y_{\ell-1},v)$, $y_{\ell+1}\in R(y_{\ell},v)$, respectively, such that
  $R(u,x_{k+1})\subsetneq R(u,x_k)$ and $R(y_{\ell+1},v)\subsetneq
  R(y_\ell,v)$. By finiteness, these sequences, eventually terminate with
  $R(u,x)=\{u,x\}$ and $R(y,v)=\{y,v\}$.
\end{proof}

\begin{corollary}
  \label{cor:pathyes}
  Let $R$ be a poset function. If $w\in R(u,v)\ne\emptyset$ then there is a
  directed path in $G_R$ from $u$ to $v$ that passes through $w$.
\end{corollary}
\begin{proof}
  Let $w\in R(u,v)$. Then Lemma~\ref{lem:upathT1} ensures that there is
  $x_1\in R(u,w)$ and $y_1\in R(w,v)$ such that $(u,x_1)\in E(G_R)$,
  $(y_1,v)$ in $E(G_R)$. Repeating the argument separately on $R(x_1,w)$
  and $R(w,y_1)$ yields edges $(x_1,x_2)\in E(R_G)$ and $(y_1,y_2)\in E(R_G)$
  and residual sets $R(x_2,w)\subseteq R(x_1,w)$ and $R(w,y_1)\subseteq
  R(w,y_2)$, which eventually yields a path
  $u,x_1,x_2\dots,x_k,w,y_h,\dots,y_2,y_1,v$ in $G_R$.
\end{proof}

\begin{lemma}
  \label{lem:pathno}
  Let $R$ be a poset function. If $R(u,v)=\emptyset$ then $G_R$ contains no
  directed path from $u$ to $v$.
\end{lemma}
\begin{proof}
  Let $R(u,v)=\emptyset$ and assume for contradiction that there is a
  directed path $P=(u,z_1,\dots,z_k,v)$.  Then $R(u,z_1)=\{u,z_1\}$ and
  $R(z_1,z_2)=\{z_1,z_2\}$. Since $R$ satisfies $\AX{tr}$ by
  Remark~\ref{remark_poset_tr}, this yield $z_1\in R(u,z_2)$ and thus
  $\{u,z_1,z_2\}\subseteq R(u,z_2)$. We proceed by induction. Assume
  $\{u,z_1,\dots,z_i\}\subseteq R(u,z_i)$. Using $R(z_i,z_{i+1}) =
  \{z_i,z_{i+1}\}$ yields $z_i\in R(u,z_{i+1})$ and by \AX{tr1}
  $R(u,z_i)\cup\{z_i,z_{i+1}\}\subseteq R(u,z_{i+1})$. Thus, in particular
  $R(u,v)$ contains all vertices of the path $P$, and thus Let
  $R(u,v)=\emptyset$ and assume, for contradiction, that there is a
  directed path $P=(u,z_1,\dots,z_k,v)$.  Then $R(u,z_1)=\{u,z_1\}$ and
  $R(z_1,z_2)=\{z_1,z_2\}$. Since $R$ satisfies $\AX{tr}$, this yield
  $z_1\in R(u,z_2)$ and thus $\{u,z_1,z_2\}\subseteq R(u,z_2)$. We proceed
  by induction. Assume $\{u,z_1,\dots,z_i\}\subseteq R(u,z_i)$. Using
  $R(z_i,z_{i+1}) = \{z_i,z_{i+1}\}$ yields $z_i\in R(u,z_{i+1})$ and by
  \AX{tr1} $R(u,z_i)\cup\{z_i,z_{i+1}\}\subseteq R(u,z_{i+1})$. Thus, in
  particular $R(u,v)$ contains all vertices of the path $P$, and thus
  $R(u,v)\ne\emptyset$, a contradiction.
\end{proof}

\begin{corollary}
  \label{cor:nopath2}
  Let $R$ be a poset function. Then the following statements are true:
  \begin{itemize}
  \item[(i)] There is a directed path connecting $u$ and $v$ in $G_R$ if
    and only if $R(u,v)\ne\emptyset$.
  \item[(ii)] If $w\notin R(u,v)$ then there is no
    directed path from $u$ to $v$ through $w$ in $G_R$.
  \item[(iii)] The graph $G_R$ is a DAG.
  \end{itemize}
\end{corollary}
\begin{proof}
  Statement (i) is an immediate consequence of Corollary~\ref{cor:pathyes}
  and Lemma~\ref{lem:pathno}. (ii) If $w\notin R(u,v)$ then the
  contrapositive of $\AX{tr}$ implies that $R(u,w)=\emptyset$ or
  $R(w,v)=\emptyset$. By Lemma~\ref{lem:pathno}, then there is no directed
  path from $u$ to $w$ or no directed path from $w$ to $v$ in $G_R$. (iii)
  following immediately from the fact that $R$ defines a partial order by
  virtue of being a poset function.
\end{proof}

Analogous to the symmetric case of the all-paths function of a undirected
graph, the all paths function of a directed graph is defined as:
\begin{definition}
  Let $G$ be a directed graph with vertex set $V$. The \emph{all-paths}
  transit function $A_G:V\times V\to 2^V$ of $G$ is defined by
  \begin{equation}
    A_G(u,v)\coloneqq
    \left\{x\in V\mid x \; \text{lies on a directed path from }   u \text{ to }v
    \right\}
  \end{equation}
\end{definition}
We are now in the position to state the main result of this section: 
\begin{theorem}
  \label{thm:poset}
  Let $R$ be a poset function, $G_R$ its underlying graph, and $A_{G_R}$ the
  all-paths transit function of $G_R$. Then $G_R$ is a DAG and $R=A_{G_R}$.
\end{theorem}
\begin{proof}
  It follows from Corollary~\ref{cor:nopath2} that $G_R$ is a DAG.
  We have $w\in A_{G_R}(u,v)$ if and only if there is a directed path from
  $u$ to $v$ through $w$ in $G_R$, which by Cor.\ref{cor:pathyes} and
  Cor.~\ref{cor:nopath2}(ii) is the case if and only $w\in R(u,v)$.  Thus
  $R(u,v)=A_{G_{R}}(u,v)$ for all $u,v\in V$.
\end{proof}
We shall return to the function $A_G$ in a more general setting in
Section~\ref{sect:daptf}.

\subsection{Shortcut free DAGs} 

The directed paths of a directed acyclic graph(DAG) $G$ define a partial order $\succeq_G$ and thus
a poset function $R_G$. It is tempting to assume that $G$ equals the
underlying graph $G_{R_G}$ of the directed transit function $R_G$. This is
not true in full generality, however, and the obstructions are
\emph{shortcuts}, i.e, transitive edges in $G$ \cite{Hellmuth:23a}. More
formally, if $P=(x_1,x_2,\dots,x_k)$ is a path of length $k\ge 2$ in $G$,
we call an edge $(x_1,x_k)$ a \emph{shortcut}. Obviously, shortcuts do not
affect the set of vertices lying on paths from $u$ to $v$, and $R_G=R_{G'}$
if $G$ and $G'$ differ only by the presence or absence of shortcut edges.
In the following we denote the set of directed graphs without shortcuts by
$\mathfrak{G}^{\circ}$.

\begin{lemma}
  \label{lem:shortcutfreeGRG}
  Let $G$ be a DAG. Then $G=G_{R_G}$ if and only if $G\in
  \mathfrak{G}^\circ$.
\end{lemma}
\begin{proof}  
  Assume that $G=G_{R_G}$. If $G\notin \mathfrak{G}^{\circ}$, then $G$
  contains a path $P=(x_1,x_2,x_3 \ldots,x_n)$ and the edge $(x_1,x_n)$.
  Then $P\subseteq R_G(x_1,x_n) \ne \{x_1,x_n\}$, and thus $(x_1,x_n)$ is
  not an edge in $G_{R_G}$. Conversely, assume that $G\neq G_{R_G}$. Note
  that $(x_1,x_2)$ is an arc in $G_{R_G}$ if and only if
  $R_G(x_1,x_2)=\{x_1,x_2\}$. That is, $(x_1,x_2)$ is an arc in $G$ .
  Therefore, $G\neq G_{R_G}$ implies that there is an edge $(x_1,x_2)$ in
  $G$ that is not in $G_{R_G}$. Thus there is $x_3\in R_{G}(x_1,x_2)$ such
  that $x_3\notin \{x_1,x_2\}$. Hence, there is a directed path from $x_1$
  to $x_2$ other than the edge $(x_1,x_2)$, and thus $(x_1,x_2)$ is
  shortcut in $G$.
\end{proof}
The graph $G_{R_G}$ is obtained from $G$ by removing all shortcuts.  For
DAG, it coincides with the unique \emph{transitive reduction} of $G$
\cite{Aho:72} and the unique minimum equivalent graph \cite{Moyles:69}.
\begin{observation}
  For every directed acyclic graph $G$ holds $G_{R_G}\in
  \mathfrak{G}^{\circ}$.
\end{observation}

\subsection{Connection to the undirected poset function}

A (symmetric) transit function $P$ has been associated with a poset
$(X,\preceq)$ by defining the order intervals $[u,v]\coloneqq \{x\in V|
u\preceq x\preceq v \text{ or } v \preceq x \preceq u\}$ and setting
$P(u,v)\coloneq [u,v]$ if $u\preceq v$ or $v\preceq u$, and if both
$u\not\preceq v$ and $v\not\preceq u$ then $P(u,v)=\{u,v\}$.  The function
$P$ has a long history and was part of folklore in the study of order
intervals and order convexity; for further references we refer to
\cite{brev2009,mathews2008,vandeVel:93}.  It satisfies a rich portfolio of
properties that admit immediate generalizations to the non-symmetric
case. In Lemma~\ref{lem:bbb}, we show that this is true in particular for
the betweenness axioms most commonly studied for symmetric transit
functions \cite{Changat:19}.
\begin{lemma}
  \label{lem:bbb}
  Let $R$ be a poset function. Then $R$ satisfies the following
  properties
  \begin{itemize}
  \item[\AX{m}] $x,y\in R(u,v)$ implies $R(x,y)\subseteq R(u,v)$
  \item[\AX{b1$_1$}] $x\in R(u,v)$ and $x\neq v$ implies $v\notin R(u,x)$
  \item[\AX{b1$_2$}] $x\in R(u,v)$ and $x\neq u$ implies $u\notin R(x,v)$
  \item[\AX{b2$_1$}] $x\in R(u,v)$ implies $R(u,x)\subseteq R(u,v)$
  \item[\AX{b2$_2$}] $x\in R(u,v)$ implies $R(x,v)\subseteq R(u,v)$      
  \item[\AX{b3$_1$}] $x\in R(u,v)$ and $y\in R(u,x)$ implies $x\in R(y,v)$
  \item[\AX{b3$_2$}] $x\in R(u,v)$ and $y\in R(x,v)$ implies $x\in R(u,y)$  
  \item[\AX{b4}]     If $x\in R(u,v)$ then $R(u,x)\cap R(x,v)=\{x\}$
  \item[\AX{Ch}] $x\in R(u,v)$ and $y\in R(x,w)$ implies $y\in R(u,w)$ or
    $y\in R(v,w)$ or $y\in R(u,v)$.
  \item[\AX{j0}] $x\in R(u,y)$ and $y\in R(x,v)$ implies $x\in R(u,v)$.
  \end{itemize}
\end{lemma}
In the following we write \AX{b1$_{1,2}$}, \AX{b2$_{1,2}$}, or \AX{b2$_2$}
if $R$ satisfies both \AX{b1$_1$} and \AX{b1$_2$}, both \AX{b2$_1$} and
\AX{b2$_2$}, or both \AX{b3$_1$} and \AX{b3$_2$}, respectively.
\begin{proof}
  If $R(x,y)=\emptyset$, there is nothing to show. Thus we assume
  $R(x,y)\ne\emptyset$. By \AX{tr2}, $x,y\in R(u,v)$ implies
  $R(u,x)\ne\emptyset$, $R(x,v)\ne\emptyset$, $R(u,y)\ne\emptyset$, and
  $R(y,v)\ne \emptyset$. By \AX{tr1}, $R(u,x)\ne\emptyset$ and $R(x,y)\ne
  \emptyset$ implies $R(x,y)\subseteq R(u,y)\ne\emptyset$, which together
  with $R(y,v)\ne \emptyset$ implies $R(u,y)\subseteq R(u,v)$, and thus
  also $R(x,y)\subseteq R(u,v)$.

  \AX{b1$_{1,2}$} is a direct consequence of \AX{q}.  The contrapositive of
  \AX{tr2} and \AX{tr1} yield $w\in R(u,v) \implies R(u,w)\cup
  R(w,v)\subseteq R(u,v)$. Together with \AX{q} implies \AX{b2$_{1,2}$}.
    
  To see \AX{b3$_1$}, we assume $x\notin R(y,v)$, whence by \AX{tr}, we
  have $R(y,x)=\emptyset$ or $R(x,v)=\emptyset$. The latter is impossible
  since $x\in R(u,v)$ and thus $R(y,x)=\emptyset$. However, we have $y\in
  R(u,x)$ and thus both $R(u,y)\ne\emptyset$ and $R(y,x)\ne\emptyset$, a
  contradiction.  \AX{b3$_2$} follows analogously.

  \AX{b4} is an immediate consequence of \AX{q}.

  $x\in R(u,v)$ and $y\in R(x,w)$ implies $R(u,x)\ne\emptyset$,
  $R(x,v)\ne\emptyset$, $R(x,y)\ne\emptyset$ and $R(y,w)\ne\emptyset$.
  Hence $y\in R(u,x)\cup R(x,y)\cup R(y,w)\subseteq R(u,w)$, and thus
  \AX{Ch} holds.
    
  Since $x\in R(u,y)$ and $y\in R(x,v)$ implies $R(u,x)\ne\emptyset$,
  $R(x,y)\ne\emptyset$ and $R(y,v)\ne\emptyset$, we obtain $x\in R(u,y)\cup
  R(y,v)\subseteq R(u,v)$ and thus \AX{j0} is satisfied.
\end{proof}

Note that the analogs of axioms \AX{b1}, \AX{b2}, and \AX{b3} come in pairs in the
directed case. We shall return to the implications of the betweenness
axioms in a more general context in Sect.~\ref{sect:geom}. Axiom \AX{Ch} first
appeared in the context of convex geometries in \cite{chva}. It may be
worth noting that one of the three alternatives appearing in this condition
is always true.  Axiom \AX{j0} was initially introduced in \cite{mcjmhm-10}
and also plays a role in the study of convex geometries.  %The importance of(Ch) and (j0).
 It remains an open question whether \AX{Ch} and \AX{j0} also play
a major role in the directed setting, given that convex geometries are
inherently undirected structures. 

\begin{lemma}
  Let $R$ be a directed transit function. Then \AX{q} implies \AX{t2a}
\end{lemma}
\begin{proof}
  Let $R(u,v)\ne\emptyset$ and $u\neq v$.  Suppose that
  $R(v,u)\ne\emptyset$. Then $\{u,v\}\subseteq R(u,v)\cap R(v,u)$. Now
  \AX{q} implies that $u=v$, a contradiction. Thus $R(v,u)=\emptyset$ and
  hence $R$ satisfies \AX{t2a}.
\end{proof}

\subsection{Rooted Directed Forests and Trees}

Rooted directed forests $\vec{F}$ and rooted directed trees $\vec{T}$ are an
interesting special case of posets. We note that the various notions of
path transit functions, which we will discuss in more detail in the
following section, coincide in $\vec{F}$ and $\vec{T}$, because any two
vertices $u$ and $v$ are connected by directed path if and only if $u$ lies
on the unique directed path from the root of a tree $\vec{T}$ (which may be
a component of a forest $\vec{F}$) to a vertex $v$, i.e., if and only if
$u$ is an ancestor of $v$.  The uniqueness of paths on $\vec{F}$ incurs
additional conditions that restrict further the directed transit
function. The main purpose of this section is to show that directed rooted
forests and trees can be characterized by suitable axioms on directed
transit functions.

Clearly the directed path transit function of $\vec{F}$, and thus
also of a tree $\vec{T}$ is a poset function. Since all paths between
vertices in $\vec{F}$ and $\vec{T}$ are unique, it is easy to see that all
transit functions defined by means of various types of paths in $\vec{F}$
and $\vec{T}$ are equal. Therefore we use the ambiguous term ''path
function'' in this section instead of referring to a specific type, such as
the all-paths function.

\begin{lemma}
  \label{lem:axiomsFor}
  The $R$ be the path function of a rooted forest $\vec{F}$.
  Then $R$ is a poset function satisfying 
  \begin{itemize}
  \item[\AX{p}] $w\in R(u,v)$ implies $R(u,w)\cup R(w,v)=R(u,v)$.
  \item[\AX{hy}] If $R(x,v)\neq \emptyset$ and $R(y,v)\neq \emptyset$ then
    $x\in R(y,v)$ or $y\in R(x,v)$.
  \end{itemize}
\end{lemma}
\begin{proof}
  \AX{p}. If $w\in R(u,v)$, then there is a unique directed path from $u$
  to $v$ in $\vec{F}$ passing through $w$.  If $w=u$ or $w=v$, or $u=v=w$
  there is nothing to show, hence assume $w\neq u \neq v$. In particular
  then, there is a unique path from $u$ to $v$, which is the concatenation
  of the unique paths from $u$ to $w$ and from $w$ to $v$, and thus
  $R(u,v)=R(u,w)\cup R(w,v)$. Hence $R$ satisfies \AX{p}.

  \AX{hy}.  Let $R(x,v)\neq\emptyset$ and $R(y,v)\neq\emptyset$, then
  clearly $x,y,v$ are in the same connected component in
  $\overrightarrow{F}$. Assume $x\notin R(y,v)$ and $y\notin R(x,v)$ then
  we have a path in $\overrightarrow{F}$ from $x$ to $v$ that does not
  include $y$, and from $y$ to $v$ that does not include $x$. Clearly, now
  either $v$ does not have a unique predecessor, or there is a third vertex
  $w\in R(x,v)\cap R(y,v)$ that does not have a unique predecessor,
  contradicting the assumption that $\overrightarrow{F}$ is a forest.
\end{proof}

\begin{observation}
  Since the path transit function $R$ of a rooted tree $\vec{T}$ in
  particular derives from a poset, we have also
  \begin{itemize}
  \item[\AX{a5*}] $R(u,x)\neq\emptyset$ and $R(x,v)\neq\emptyset$ implies
    $R(u,x)\cup R(x,v) = R(u,v)$.
  \item[\AX{a5}]  $R(u,x)\cap R(x,v)=\{x\}$ implies
    $R(u,x)\cup R(x,v)=R(u,v)$
  \end{itemize}
\end{observation}
\begin{proof}
  To see this, note that $R(u,x)\ne\emptyset$ and $R(x,v)\ne\emptyset$ by
  \AX{tr1} implies $x\in R(x,v)\cap R(x,v)$ and hence \AX{p} and \AX{tr2}
  yields \AX{a5*}; moreover, \AX{a5*} implies \AX{a5}.
\end{proof}
Axiom \AX{a5} appears, e.g.,\ in \cite{Changat:01}. We shall see below that
axiom \AX{p} yields unique paths for the underlying graph of a poset
function.

\begin{lemma}
  \label{lem:upathT-T}
  Let $R$ be a poset function satisfying \AX{p}.  If $R(u,v)\ne\emptyset$
  and $u\ne v$ then there are unique $x,y\in R(u,v)$ such that
  $R(u,x)=\{u,x\}$ and $R(y,v)=\{y,v\}$.
\end{lemma}
\begin{proof}
  The existence of $x$ and $y$ has been established in
  Lemma~\ref{lem:upathT1} for all poset functions. It remains to show
  uniqueness. To this end, consider $x,x'\in R(u,v)$ such that
  $R(u,x)=\{u,x\}$ and $R(u,x')=\{u,x'\}$. Assume that $x'\ne x$. Then
  \AX{p} implies $x'\in R(x,v) = R(x,x')\cup R(x',v)$. Since we also have
  $R(u,v)=R(u,x')\cup R(x',v)$ by \AX{p} and $R(u,x')\cap R(x',v)=\{x'\}$
  by \AX{q}(by Lemma~\ref{lemma-axiom-q}), we conclude that $x'\in R(u,x)$,
  a contradiction. Thus $x=x'$. An analogous argument shows that $y,y'\in
  R(u,v)$ with $R(y,v)=\{y,v\}$ and $R(y',v)=\{y',v\}$ implies
  $y=y'$.\\ Consider $x,x'\in R(u,v)$ such that $R(u,x)=\{u,x\}$ and
  $R(u,x')=\{u,x'\}$. Assume that $x'\ne x$. Then \AX{p} implies $x'\in
  R(x,v) = R(x,x')\cup R(x',v)$. Moreover, $x\in R(u,v)\implies R(u,x)\cup
  R(x,v)=R(u,v)$ and $x'\in R(x,v)\implies R(x,x')\cup R(x'v)=R(x,v)$.  By
  similar argument we can show that $x\in R(x',v)\implies R(x',x)\cup
  R(x,v)=R(x'v)$.  Now either $R(x,x')=\emptyset$ or $R(x',x)=\emptyset$.
  W.l.o.g. let $R(x,x') =\emptyset$ and so $R(x',x)\neq\emptyset \implies
  x,x'\in R(x',x)$. This in turn implies $R(x',x)\cap R(x,v)=\{x,x'\}$, a
  contradiction to the axiom \AX{q}.  Hence $x=x'$.  In a similar manner it
  can be shown that $y=y'$.
\end{proof}

An immediate consequence of Lemma~\ref{lem:upathT-T} is the existence of
unique paths in the underlying graph $G_R$ of a poset function $R$ that
satisfies \AX{p}:
\begin{corollary}
  \label{cor:pathT}
  Let $R$ be a poset function satisfying \AX{p}. Then for any two
  $u,v\in V$ for which $R(u,v)\neq\emptyset$, $R(u,v)$ consists of
  precisely the vertices that lie on the unique path from $u$ to $v$ in
  $G_R$.
\end{corollary}

A \emph{mangrove} or \emph{strongly unambigious graphs} \cite{Allender:98},
also termed \emph{multitree} \cite{Furnas:94} is a DAG $G$ with at most one
directed path between any two vertices. Equivalently, the reachability
relation $\succeq$ of such a DAG $G$ is diamond-free poset, i.e., there are
no four distinct elements $a,b,c,d\in V$ such that $a\succeq b \succeq d$,
$a\succeq c \succeq d$, and $b\parallel c$ \cite{Furnas:94}. As an
immediate consequence of Theorem~\ref{thm:poset} and Cor.~\ref{cor:pathT}
we obtain:
\begin{corollary}
  Let $R$ be a poset function that satisfies \AX{p}, $G_R$ its underlying
  graph, and $R_{G_R}$ the path transit function of $G_R$. Then $G_R$ is a
  mangrove and $R=R_{G_R}$.
\end{corollary}

The intuition behind condition \AX{hy} is to forbid ``hybrid vertices'',
i.e., the confluence of multiple paths in a vertex $v$.
\begin{lemma}
  \label{lem:degreeT}
  Let $R$ be a poset function satisfying \AX{p} and \AX{hy}. Then all
  vertices in $G_R$ have in-degree at most $1$.
\end{lemma}
\begin{proof}
  Let $R$ be a poset function that satisfies \AX{p} and \AX{hy}. Assume
  $G_R$ contains a vertex $v$ with in-degree larger than one. Then we have
  two vertices $x,y$ such that $R(x,v)=\{x,v\}$ and $R(y,v)=\{y,v\}$ but,
  by \AX{hy}, either $x\in R(y,v)$ or $y\in R(x,v)$, a contradiction.
\end{proof}

As an immediate consequence of Corollary~\ref{cor:pathT} and
Lemma~\ref{lem:degreeT} we obtain
\begin{corollary}
  \label{cor:forestroot}
  Let $R$ be a transit function that satisfies \AX{p} and \AX{hy}. Then
  every weakly connected component of $G_R$ contains a unique source
  vertex.
\end{corollary}
The unique source vertex is the root of the trees that form the weakly
connected components of $\vec{F}$.

\begin{theorem}
  Let $R$ be a poset function that satisfies \AX{p} and \AX{hy}, $G_R$ its
  underlying graph, and $R_{G_R}$ the path transit function of $G_R$. Then
  $G_R$ is a directed rooted forest and $R_{G_R}=R$.
\end{theorem}
\begin{proof}
  By Theorem~\ref{thm:poset}, $G_R$ is a DAG, by Corollary~\ref{cor:pathT}
  we know that any path between two vertices in $G_R$ is unique, and by
  Corollary~\ref{cor:forestroot} we know that every weakly connected
  component in $G_R$ contains a unique root. Since $R$ is a poset function,
  Theorem~\ref{thm:poset} implies $R_{G_R}=R$.
\end{proof}
  
Let us now turn to trees. Here, we have to establish a unique root and
connectedness.
\begin{lemma}
  \label{lem:axiomsT}
  Let $R$ be the path function of a rooted tree $\vec{T}$.  Then $R$ is a
  poset function satisfying \AX{p} and
  \begin{itemize}
  \item[\AX{r}] $R(u,v)=R(v,u)=\emptyset$ implies that there is $w\in V$
    such that $R(w,u)\neq \emptyset$ and $R(w,v)\neq\emptyset$.
  \end{itemize}
\end{lemma}
\begin{proof}
  Since every rooted tree is in particular a rooted forest,
  Lemma~\ref{lem:axiomsFor} implies that \AX{p} holds. Suppose
  $R(u,v)=R(v,u)=\emptyset$ with $u,v\in V(\vec{T})$, i.e., there is
  neither a path from $u$ to $v$ or from $v$ to $u$ in $\vec{T}$. Since
  $\vec{T}$ is a rooted tree there is a root vertex $w\in V(\vec{T})$ such
  that there is a path from $w$ to $u$ and $v$. Consequently
  $R(w,u)\neq\emptyset$ and $R(w,v)\neq\emptyset$ and \AX{r} holds.
\end{proof}
Intriguingly, a poset function $R$ that satisfies \AX{p} together with
\AX{r} yields a graph $G_R$ that does not contain any hybrid vertices.

\begin{lemma}
  Let $R$ be a poset function that satisfies \AX{p} and \AX{r}. Then $R$
  satisfies \AX{hy}.
\end{lemma}
\begin{proof}
  Let $x,y,v\in V$ such that $R(x,v)\neq\emptyset$ and
  $R(y,v)\neq\emptyset$. Without loss of generality assume
  $R(x,y)\neq\emptyset$, then \AX{tr} implies $y\in R(x,v)$.
		
  Now consider the case that both $R(x,y)=\emptyset$ and $R(y,x)=\emptyset$
  and hence $x\ne y$. By \AX{r}, there is $w\in V$ such that
  $R(w,x)\neq\emptyset$ and $R(w,y)\neq\emptyset$, and hence \AX{tr} implies
  $x,y\in R(w,v)$. By \AX{p} we have $R(w,v)=R(w,x)\cup
  R(x,v)=R(w,y)\cup R(y,v)$. If $x\in R(y,v)$ or $y\in R(x,v)$ there is
  nothing to show. Hence assume $y\in R(w,x)$ and $x\in R(w,y)$.  Using
  \AX{p} again we obtain $R(w,x)=R(w,y)\cup R(y,x)$. Moreover, by \AX{q}
  $R(w,x)=R(w,y)\cap R(y,x)=\{y\}$ and thus $x\notin R(w,y)$, a
  contradiction. 		
\end{proof}

\begin{lemma}
  \label{lem:rootT}
  Let $R$ be a poset function that satisfies \AX{r}, then $G_R$ contains
  a unique source vertex.
\end{lemma}
\begin{proof}
  Let $R$ be a poset function on $V$ that satisfies \AX{r}. If for every
  $u,v$ we have $u=v$, or $R(u,v)\ne\emptyset$, or $R(v,u)\ne \emptyset$,
  then $\succeq_R$ is a total order, and thus there is a unique $w\in V$
  such that $R(w,u)\ne\emptyset$ for all $u\in V$, i.e., there is unique
  source vertex. Now suppose, for contradiction, that no such vertex
  exists. Then $\succeq_R$ cannot be a total order, and thus there is a
  pair of vertices $u,v\in V$ such that $R(u,v)=R(v,u)=\emptyset$. Let
  $U_1=\{u,v\}$. By \AX{r}, there is $u_1\in V$ with $R(u_1,u)\ne\emptyset$
  and $R(u_1,v)\ne\emptyset$, and thus $u_1\notin U_1$. By assumption,
  there is $v_1\in V$ such that $R(u_1,v_1)=\emptyset$ and
  $R(v_1,u_1)=\emptyset$, which in particular implies $v_1\notin U_1$. By
  \AX{r}, there is $u_2\in V$ such that $R(u_2,u_1)\ne\emptyset$ and
  $R(u_2,v_1)\ne\emptyset$, thus by \AX{t0} also $R(u_2,x)\ne\emptyset$ for
  all $x\in U_1\cup\{u_1,v_1\}\eqqcolon U_2$. Repeating this argument
  yields a sequence of vertices $u_i,v_i\in V\setminus U_{i-1}$ such that
  $R(u_i,v_i)=R(v_i,u_i)=\emptyset$ and $R(u_{i},x)\ne\emptyset$ for all
  $x\in U_{i-1}$ and $U_i=U_{i-1}\cup\{u_{i-1},v_{i-1}\}$. In particular
  the set $U_i$ grows in each step, which is impossible since $V$ is
  finite. Hence, there must be a source vertex.
\end{proof}
The existence of a unique source vertex in a graph in particular implies
that this graph is (at least) weakly connected.

\begin{lemma}
  \label{lem:connectedT}
  Let $R$ be a poset function that satisfies \AX{r} and \AX{p}. Then $G_R$
  is weakly connected.
\end{lemma}
\begin{proof}
  Let $u,v\in V$ such that $R(u,v)=R(v,u)=\emptyset$ then we know by \AX{r}
  that there is a $w\in V$ such that $R(w,u)\neq\emptyset$ and
  $R(w,v)\neq\emptyset$. By Corollary~\ref{cor:pathT}, there is a path from
  $w$ to $u$ in $G_R$ and from $w$ to $v$ in $G_R$. Hence, $w,u,$ and $v$
  are contained within the same (weakly) connected component in $G_R$.
		
  If $R(u,v)\neq\emptyset$ and $R(v,u)=\emptyset$, we use
  Corollary~\ref{cor:pathT} again to show that there is a path from $u$ to
  $v$ in $G_R$ (analogously for $R(v,u)\neq\emptyset$ and
  $R(u,v)=\emptyset$). Since trivially \AX{t2a} prohibits symmetric edges,
  the case that $R(u,v)\neq\emptyset$ and $R(v,u)\neq\emptyset$ cannot
  occur. Summarizing, for any two vertices in $G_R$, there is either a path
  between them, or there is a third vertex $w$ such that there is a path
  from $w$ to $u$ and from $w$ to $v$. Consequently, $G_R$ is weakly
  connected.
\end{proof}

So far, we have established that the underlying graph $G_R$ of a poset
function $R$ that satisfies \AX{p} and \AX{r} is weakly connected, contains
a unique source vertex, and there is at most one directed path between any
two vertices of $G_R$. This puts us into position to state the main result
of this subsection.
\begin{theorem}
  Let $R$ be a poset function that satisfies \AX{r} and \AX{p}, $G_R$ its
  underlying graph, and $R_{G_R}$ the all-path transit function of
  $G_R$. Then $G_R$ is a directed rooted tree and $R=R_{G_R}$.
\end{theorem}
\begin{proof}
  By Theorem~\ref{thm:poset}, $G_R$ is a DAG and by
  Lemma~\ref{lem:connectedT} $G_R$ is weakly connected. We know by
  Corollary~\ref{cor:pathT} that every path from one vertex to another in
  $G_R$ is unique. The existence of a root in $G_R$
  (by Lemma~\ref{lem:rootT}), from which there exists a unique path to every
  other vertex in $V$ together with the fact that $G_R$ is a DAG implies
  that the root is unique. Consequently, $G_R$ is a rooted tree.  Since $R$
  is a poset function, Theorem~\ref{thm:poset} implies $R_{G_R}=R$.
\end{proof}

\section{Directed transit functions of general graphs}
\label{sect:reach}

\subsection{Reachability}

Reachability in its most general form can be understood as a reflexive,
transitive relation $\to$ on a set $V$. Naturally, we set
\begin{equation}
  R_{\to}(u,v) = \begin{cases}
    \{x\in V| u \to x \text{ and } x \to v \} & \text{if } u\ne v \\
    \{u\}                                     & \text{if } u = v \\
    \end{cases}
\end{equation}
Note that $R_{\to}(u,v)=\emptyset$ if and only if $u\not\to v$, i.e., if $v$
is not reachable from $u$. One easily checks that $R_{\to}$ satisfies
\AX{t1}, \AX{t3}, \AX{tr1}, and \AX{tr2}, using the same arguments as in
Lemma~\ref{lem:poset1}.

The statements of \AX{tr1} and \AX{tr2} do not extend to the case $u=v$ in
the absence of the anti-symmetry axiom \AX{t2a}, which is not assumed
here. Therefore, \AX{q} no longer holds in this more general setting.

\begin{figure}
  \begin{center}
    \begin{tikzpicture}
      \tikzset{vertex/.style = {shape=circle,draw,minimum size=1.0em}}
      \tikzset{edge/.style = {->,> = latex'}}
      \node[vertex] (1) at  (0,0) {$u$};
      \node[vertex] (2) at  (1.25,2) {$v$};
      \node[vertex] (3) at  (2.5,0) {$w$};
      
      \draw[edge] (1) to [bend left=15] (2);
      \draw[edge] (2) to [bend left=15] (1);
      
      \draw[edge] (2) to [bend left=15] (3);
      \draw[edge] (3) to [bend left=15] (2);
      
      \draw[edge] (1) to [bend left=15] (3);
      \draw[edge] (3) to [bend left=15] (1);
    \end{tikzpicture}
  \end{center}
  \caption{The reachability function of this digraph satisfies
      $R_{\to}(x_1,x_2)=\{u,v,w\}$ for any two distinct vertices
      $x_1,x_2\in\{u,v,w\}$. Thus $G_{R_{\to}}$ is edge-less and hence
      $G\not\cong G_{R_{\to}}$.}
  \label{fig:G!=G.to}
\end{figure}
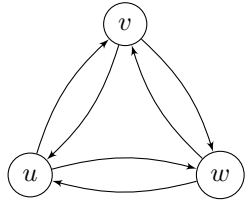

Let $G=(V,E)$ be a connected graph. Then for every $u$ and every $v\ne u$,
we have $u\to x$ and $x\to v$, and thus $x\in R_{\to}(u,v)$. That is,
$R_{\to}(u,v)=V$ for $u\ne v$ and $R_{\to}(u,u)=\{u\}$. The example in
Fig.~\ref{fig:G!=G.to} show that in general the underlying graph of the
reachability relation does not recover the original graph.

Reachability is closely related to notions of connectedness.  We see that
$V$ is \emph{strongly connected} if $R(u,v)\ne\emptyset$ for all $u,v\in V$
and \emph{unilaterally connected} if $R(u,v)\ne\emptyset$ or
$R(v,u)\ne\emptyset$ holds for all $u,v\in V$.
  
\begin{lemma}
  Let $R$ be a directed transit function on $V$. Then the binary relation
  $\sim$ on $V$ defined by $u\sim v$ if $R(u,v)\ne\emptyset$ and
  $R(v,u)\ne\emptyset$ is an equivalence relation.
\end{lemma}
\begin{proof}
  Symmetry is inherent in the definition and by \AX{t3}, $\sim$ is
  reflexive.  Suppose $u\sim w$ and $w\sim v$. Then $R(u,w)\ne\emptyset$
  and $R(w,v)\ne\emptyset$, which by \AX{t0} implies $R(u,v)\ne\emptyset$;
  analogously, $R(v,w)\ne\emptyset$ and $R(w,u)\ne\emptyset$ by \AX{t0}
  implies $R(v,u)\ne\emptyset$, and thus $v\sim u$, i.e., $\sim$ is an
  equivalence relation.
\end{proof}
The equivalence classes of $\sim$ can be viewed as strongly connected
components w.r.t.\ the directed transit function $R$. The set $V$ is
partitioned into a set $\mathcal{C}$ of strongly connected components if
and and only if $R(u,v)\ne\emptyset$ implies $R(v,u)\ne\emptyset$ for all
$u,v\in C$ with $C\in\mathcal{C}$. Moreover, $V$ is strongly connected
w.r.t.\ $R$ if and only if $R(u,v)\ne\emptyset$ for all $u,v\in V$.

\begin{lemma}
  \label{lem:nopath}
  Let $R$ be a directed transit function. Then $R(u,v)=\emptyset$ implies
  that there is no path from $u$ to $v$ in the underlying graph $G_R$.
\end{lemma}
\begin{proof}
  Let $C_+(u)\coloneqq\{w\in V| R(u,w)\ne\emptyset\}$ and
  $C_-(v)\coloneqq\{w\in V| R(w,v)\ne\emptyset$\}. By \AX{t0},
  $R(u,v)=\emptyset$ implies $C_+(u)\cap C_-(v)=\emptyset$. Suppose there
  is an edge $(p,q)\in E(G_R)$ with $p\in C_+(u)$ and $q\in C_-(u)$. Then
  $R(p,q)=\{p,q\}\ne\emptyset$ and thus $R(u,p)\ne\emptyset$ and
  $R(p,q)\ne\emptyset$ implies $R(u,q)\ne\emptyset$, i.e., $q\in C_+(u)$, a
  contradiction. Thus $G_R$ cannot contain an edge connecting $u$ with any
  vertex $q\notin C_+(u)$, and thus there cannot be a path in $G$ for $u$
  to any $q\notin C_+(u)$.
\end{proof}
Lemma~\ref{lem:nopath} extends Lemma~\ref{lem:pathno} to all directed
transit functions. The converse is not true in general for directed transit
functions since $R$ need not have any transit sets of size $2$.

\subsection{Geometric directed transit functions} 
\label{sect:geom}

It is interesting to note that key implications between properties of poset
functions remain true when the condition of the anti-symmetry axiom
\AX{t2a} is lifted.  In particular, it is well-known for the case of symmetric transit function
that (b3) implies (b4) implies (b1). This statement generalizes to the
directed transit functions:
\begin{lemma}\label{lemma_b3_implies_b1}
  Let $R$ be a directed transit function satisfying \AX{tr2}. If $R$
  satisfies \AX{b3$_{1}$} and \AX{b3$_{2}$}, then it also satisfies
  \AX{b4}, \AX{b1$_{1}$}, and \AX{b1$_{2}$}.
\end{lemma}
\begin{proof}      
  Suppose $x\in R(u,v)$. By \AX{tr2}, $R(u,x)\ne\emptyset$ and
  $R(x,v)\ne\emptyset$ and hence $x\in R(u,x)\cap R(x,v)$. Let $y\in
  R(u,x)$ and $y\in R(x,v)$. Then by \AX{b3$_1$} $x\in R(y,v)$ and by
  \AX{b3$_2$} $x\in R(u,y)$.  Applying \AX{b3$_1$} to $x\in R(u,y)$ and
  $y\in R(u,x)$ with $v=y$ implies $x\in R(y,y)=\{y\}$ and thus $x=y$, and
  thus \AX{b4} holds.  Now consider $x\in R(u,v)\ne\{u,v\}$. Then
  $R(u,x)\cap R(x,v)=\{x\}$ by \AX{b4} and thus $u\in R(u,x)$ implies
  $u\notin R(x,v)$, and $v\in R(v,x)$ implies $v\notin R(u,x)$, i.e.,
  \AX{b1$_{1}$} and \AX{b1$_{2}$} hold.
\end{proof}

For symmetric transit functions without empty transit sets, the notion of a
geometric (transit) functions was introduced in \cite{Nebesky:01} as those
that satisfy axioms \AX{b2} and \AX{b3}. Since then, the term has been in
common use; see also \cite{Bandelt:96,vandeVel:93}.  In the following we
combine \AX{b2$_{1}$} and \AX{b2$_{2}$} defined in Lemma~\ref{lem:bbb} into
a single axiom
\begin{itemize}
\item[\AX{b2}] $w\in R(u,v)$ implies $R(u,w)\cup R(w,v)\subseteq R(u,v)$
\end{itemize}
Recalling \AX{tr0} from Observation~\ref{obs:tr0}, it is easy to check
\begin{observation}
  A directed transit function $R$ satisfies \AX{tr1} if and only it satisfies
  \AX{tr0} and \AX{b2}.
\end{observation}

Axiom \AX{tr0} is clearly true for the reachability transit function, and
for the directed all-paths transit function in the setting of partial
orders. In general it fails for path-based transit functions since the
concatenation of a $u,w$-path and a $w,v$-path is, in general, only
guaranteed to contain a $u,v$-path. However, this $u,v$-path need not pass
through the concatenation point $w$. In \cite{Changat:17}, therefore the
weaker condition \AX{t0} appears to ensure consistency of connected
components in undirected graphs. Axiom \AX{t0} is an immediate
consequence of \AX{tr1}. Notably, even in conjunction with \AX{tr2} and
\AX{b3$_{1,2}$}, \AX{t0} does not imply \AX{tr1}.

A natural generalization of Nebesk{\'y}'s notion of a geometric transit
function for (not necessarily connected) structures is the following:
\begin{definition}
  A directed transit function $R$ is \emph{geometric} if it satisfies
  \AX{tr2}, \AX{b2} and \AX{b3$_{1,2}$}. Moreover, we say that directed
  transit function $R$ is \emph{weakly geometric} if it satisfies \AX{tr2},
  \AX{b2} and \AX{b1$_{1,2}$}.
\end{definition}

It follows from Lemma~\ref{lemma_b3_implies_b1} that geometric directed
transit function is always weakly geometric. The converse is not true in
general, as shown by examples \ref{ex-b1_not_b31} and \ref{ex-b1_not_b32}. 
\begin{example}\label{ex-b1_not_b31}
  \AX{b1$_{1,2}$}, \AX{b2}, but not \AX{b3$_{1}$}.\newline
  $V=\{a,b,c,d\}$ and $R(a,d)=\{a,b,c,d\}$, $R(a,c)=\{a,b,c\}$,
  $R(a,b)=\{a,b\}$, $R(b,c)=\{b,c\}$, $R(c,d)=\{c,d\}$, $R(b,d)=\{b,d \}$,
  $R(x,x)=\{x\}$ for all $x\in V$, and $R(u,v)=\emptyset$ otherwise.
\end{example}
\begin{example}\label{ex-b1_not_b32}
  \AX{b1$_{1,2}$}, \AX{b2}, but not \AX{b3$_2$}.\newline
  $V=\{a,b,c,d\}$ and $R(a,d)=\{a,b,c,d\}$, $R(b,d)=\{b,c,d\}$,
  $R(a,b)=\{a,b\}$, $R(b,c)=\{b,c\}$, $R(c,d)=\{c,d\}$, $R(a,c)=\{a,c\}$,
  $R(x,x)=\{x\}$ for all $x\in V$, and $R(u,v)=\emptyset$ otherwise.
\end{example}

The following examples of functions $R:V\times V\to 2^V$ satisfying \AX{t1}
and \AX{t3} show that the axioms \AX{t0}, \AX{b1$_1$}, \AX{b1$_2$},
\AX{b2}, and \AX{tr2} are independent of each other, i.e., weakly geometric
transit functions are not defined by any subset of this collection of
axioms.
\begin{example}\label{ex:but tr2}
  \AX{b1$_1$}, \AX{b1$_2$}, \AX{b2}, and \AX{t0} but not \AX{tr2}. \newline
  $V=\{a,b,c\}$ and $R(a,b)=\{a,c,b\}$, $R(a,c)=\{a,c\}$,
  $R(b,c)=\{b,c\}$, $R(x,x)=\{x\}$ for all $x\in V$,
  and $R(u,v)=\emptyset$, otherwise.
\end{example}
\begin{example}\label{ex:but t0}
  \AX{b1$_1$}, \AX{b1$_2$}, \AX{b2}, and \AX{tr2} but not \AX{t0}.\newline
  $V=\{a,b,c\}$ and $R(a,b)=\{a,c,b\}$, $R(b,a)=\emptyset$,
  $R(x,x)=\{x\}$ for all $x\in V$, and $R(u,v)=\{u,v\}$ otherwise.
\end{example}
\begin{example}\label{ex:but b1}
  \AX{b1$_2$}, \AX{b2}, \AX{t0}, and \AX{tr2} but not \AX{b1$_1$}.\newline
  $V=\{a,b,c,d\}$ and $R(a,b)=\{a,c,d,b\}$, $R(a,c)=\{a,b,d,c\}$,
  $R(a,d)=\{a,d\}$, $R(b,c)=\{b,c\}$, $R(c,b)=\{c,b\}$, $R(d,b)=\{d,b\}$,
  $R(d,c)=\{d,c\}$, $R(x,x)=\{x\}$ for all $x\in V$, and $R(u,v)=\emptyset$
  otherwise.
\end{example}
\begin{example}\label{ex:but b1_2}
  \AX{b1$_1$}, \AX{b2}, \AX{t0}, and \AX{tr2} but not \AX{b1$_2$}.\newline
  Let $V=\{a,b,c,d\}$ and $R(a,b)=\{a,b,c\}=R(c,b)$, $R(a,c)=\{a,c\}$,
  $R(a,d)=\{a,d\}$, $R(b,c)=\{b,c\}$, $R(b,d)=\{b,d\}$, $R(c,d)=\{c,d\}$,
  $R(x,x)=\{x\}$ and $R(u,v)=\emptyset$ otherwise.
\end{example}
\begin{example}\label{ex:but b2}
  \AX{b1$_1$}, \AX{b1$_2$}, \AX{t0}, and \AX{tr2} but not \AX{b2}.\newline
  $V=\{a,b,c,d\}$ and $R(a,b)=\{a,c,b\}$, $R(a,c)=\{a,d,c\}$,
  $R(a,d)=\{a,d\}$, $R(c,b)=\{c,b\}$, $R(d,b)=\{d,b\}$, $R(d,c)=\{d,c\}$,
  $R(x,x)=\{x\}$ for all $x\in V$ and $R(u,v)=\emptyset$ otherwise.
\end{example}

Similarly, the axioms in the definition of geometric directed transit
functions are independent:
\begin{example} \AX{t0}, \AX{b3$_{1,2}$}, \AX{b2} but not \AX{tr2}.
  \newline
  $V=\{a,b,c,d\}$ and $R(a,d)=\{a,b,c,d\}$, $R(a,c)=\{a,b,c\}$,
  $R(a,b)=\{a,b\}$, $R(b,c)=\{b,c\}$, $R(c,d)=\{c,d\}$,
  $R(b,d)=\{b,c,d\}$, $R(d,c)=\{d,c\}$, $R(d,b)=\{d,c,b\}$,
  $R(x,x)=\{x\}$ for all $x\in V$, and $R(u,v)=\emptyset$ otherwise.
\end{example}
\begin{example}
  \AX{tr2}, \AX{b3$_{1,2}$}, \AX{b2} but not \AX{t0}.
  \newline
  $V=\{a,b,c,d\}$ and $R(a,d)=\{a,b,c,d\}$, $R(a,c)=\{a,b,c\}$,
  $R(a,b)=\{a,b\}$, $R(b,c)=\{b,c\}$, $R(c,d)=\{c,d\}$, $R(b,d)=\{b,c,d\}$,
  $R(d,c)=\{d,c\}$, , $R(c,b)=\{c,b\}$, $R(x,x)=\{x\}$ for
  all $x\in V$, and $R(u,v)=\emptyset$ otherwise.
\end{example}
\begin{example}
  \AX{tr2}, \AX{b3$_{1,2}$}, \AX{t0}, \AX{b2$_2$} but not \AX{b2$_1$}.
  \newline
  $V=\{a,b,c,d\}$ and $R(a,d)=\{a,b,c,d\}$, $R(a,c)=\{a,b,c\}$,
  $R(a,b)=\{a,b\}$, $R(b,c)=\{b,c\}$, $R(c,d)=\{c,d\}$, $R(b,d)=\{b,c,d\}$,
  $R(d,b)=\{d,c,b\}$, $R(d,c)=\{a,d,c\}$, $R(x,x)=\{x\}$ for all $x\in V$,
  and $R(u,v)=\emptyset$ otherwise.
\end{example}
\begin{example}
  \AX{tr2}, \AX{b3$_{1,2}$}, \AX{t0}, \AX{b2$_1$} but not \AX{b2$_2$}.
  \newline  
  $V=\{a,b,c,d\}$ and $R(a,d)=\{a,b,c,d\}$, $R(a,c)=\{a,b,c\}$,
  $R(a,b)=\{a,b\}$, $R(b,c)=\{b,c\}$, $R(c,d)=\{c,d\}$, $R(b,d)=\{b,c,d\}$,
  $R(d,b)=\{d,c,b\}$, $R(c,b)=\{a,c,b\}$, $R(x,x)=\{x\}$ for all $x\in V$,
  and $R(u,v)=\emptyset$ otherwise.
\end{example}
\begin{example}
  \AX{t0}, \AX{b3$_{2}$}, \AX{b2}, \AX{tr2} but not 
  \AX{b3$_{1}$} 
  \newline
  $X=\{a,b,c,d\}$ and $R(a,b)=V$, $R(a,c)=\{a,c,d\}$, $R(a,d)=\{a,c,d\}$,
  $R(b,c)=\{b,c\}$, $R(c,d)=\{c,d\}$, $R(d,c)=\{d,c\}$, $R(c,b)=\{c,b\}$,
  $R(x,x)=\{x\}$ for all $x\in V$, and $R(u,v)=\emptyset$ otherwise.
\end{example}
\begin{example}
  \AX{t0}, \AX{b3$_{1}$}, \AX{b2}, \AX{tr2} but not \AX{b3$_{2}$}.
  \newline
  $X=\{a,b,c,d\}$ and $R(a,b)=V$, $R(c,b)=\{b,c,d\}$, $R(a,d)=\{a,d\}$,
  $R(b,c)=\{b,c\}$, $R(c,d)=\{c,d\}$, $R(d,c)=\{d,c\}$, $R(c,b)=\{c,b\}$,
  $R(x,x)=\{x\}$ for all $x\in V$, and $R(u,v)=\emptyset$ otherwise
\end{example}

The reachability and poset functions discussed above are special cases of
geometric transit functions:
\begin{observation}
  Let $R$ be a geometric transit function. Then:
  \begin{itemize}
  \item[(i)] $R$ is a reachability transit function if and only if it
    satisfies \AX{tr0}.
  \item[(ii)] $R$ is a poset function if and only if it satisfies
    \AX{tr0} and \AX{t2a}.
  \end{itemize}
\end{observation}

\begin{theorem}\label{thm:connectedness}
  Let $R$ be a weakly geometric directed transit function.  Then,
  $R(u,v)\ne\emptyset$ if and only if there exists a directed path from $u$
  to $v$ in $G_R$ that is contained in $R(u,v)$. If $R$ in addition
  satisfies \AX{b4}, then for every $w\in R(u,v)$ such a path can be found
  that runs through $w$. 
\end{theorem}
\begin{proof}
  If $R(u,v)=\emptyset$, then Lemma~\ref{lem:nopath} rules out the
  existence of path from $u$ to $v$ in $G_R$. In the following we assume
  $R(u,v)\neq \emptyset$. If $u=v$ then $R(u,u)=\{u\}$ and there is nothing
  to show. Thus we assume $u\neq v$ and proceed by induction over
  $|R(u,v)|$ that there is a directed path from $u$ to $v$ in $G_R$. The
  base case $|R(u,v)|=2$ is trivial since we have $R(u,v)=\{u,v\}$ and thus
  $(u,v)\in E(G_R)$. For $|R(u,v)|=3$ let $R(u,v)=\{u,v,x\}$.  We have
  $R(u,x)\subseteq R(u,v)$ and $R(x,v)\subseteq R(u,v)$ by \AX{b2} and
  $R(u,x)\subsetneq R(u,v)$ and $R(x,v)\subsetneq R(u,v)$ by \AX{b1$_1$}
  and \AX{b1$_2$}, respectively, while \AX{tr2} implies that both
  $R(u,x)\ne\emptyset$ and $R(x,v)\ne\emptyset$. Therefore $R(u,x)=\{u,x\}$
  and $R(x,v)=\{x,v\}$, and thus $(u,x,v)$ is path in $G_R$.

  Now, assume that assertion holds for all $u'$ to $v'$ such that
  $2\le|R(u',v')|< k$ and assume $|R(u,v)|=k$. Thus there is $w\in
  R(u,v)\setminus\{u,v\}$ and hence $R(u,w)\ne\emptyset$ and
  $R(w,v)\ne\emptyset$. From \AX{b2} and \AX{b1$_1$} we obtain
  $R(u,w)\subsetneq R(u,v)$ since $v\notin R(u,w)$, and thus $|R(u,w)|<k$.
  Analogously, \AX{b2} and \AX{b1$_2$} yields $R(w,v)\subsetneq R(u,v)$ and
  $|R(w,v)|<k$.  By induction hypothesis, therefore, there is a directed
  path from $u$ to $w$ contained in $R(u,w)$ and a directed path from $w$
  to $v$ contained in $R(w,v)$. The concatenation of these two directed
  paths contains a path from $u$ to $v$ that, by \AX{b2}, is contained in
  $R(u,v)$.

  If $R$ satisfies \AX{b4}, then $R(u,w)\cap R(w,v)=\{w\}$ and thus the
  concatenation of any directed path from $u$ to $w$ in $R(u,w)$ with any
  directed path from $w$ to $v$ in $R(w,v)$ is again a directed path. Thus
  there is directed path from $u$ to $v$ in $R(u,v)$ that runs through
  $w$.
\end{proof}
Theorem.~\ref{thm:connectedness} extends Cor.\ref{cor:pathyes} from poset
functions to geometric transit functions. The assumption of axiom \AX{b4}
in the second part of Theorem~\ref{thm:connectedness} is essential:
\begin{example}
  Let $V=\{u,x,w,v\}$ and $R(u,v)=\{u,x,w,v\}$, $R(u,w)=\{u,x,w\}$,
  $R(w,v)=\{w,x,v\}$, $R(u,x)=\{u,x\}$, $R(x,w)=\{x,w\}$, $R(w,x)=\{w,x\}$,
  $R(x,v)=\{x,v\}$, $R(y,y)=\{y\}$ for all $y\in V$, and $R(y,z)=\emptyset$
  otherwise. The directed transit function $R$ satisfies \AX{b1$_1$},
  \AX{b1$_2$}, \AX{b2}. However \AX{b4}, and hence \AX{b3$_1$} and
  \AX{b3$_2$} do not hold. We have $w\in R(u,v)$, there is a path from $u$
  to $v$ in $R(u,v)$, namely $(u,x,v)$, but there is no path from $u$ to
  $v$ passing through $w$.
\end{example}

\subsection{Directed all-path transit function}
\label{sect:daptf}

Here we continue with more properties of the all-paths functions $A_G$ on
an arbitrary graph $G$. 
\begin{lemma}
  \label{lem:propdir}
  Let $G$ be a directed graph. Then $R=A_G$ is a directed transit function
  satisfying \AX{tr2} as well as 
  \begin{itemize}
  \item[\AX{b5}] $R(u,w)\cap R(w,v)=\{w\}$ implies $R(u,w)\cup R(w,v)
    \subseteq R(u,v)$.
  \item[\AX{Mod}] $R(u,w)\neq\emptyset$ and $R(w,v)\neq\emptyset$
    %and $R(u,w)\neq\emptyset$
    implies $R(u,v)\cap R(u,w)\cap R(w,v) \neq \emptyset$.
  \end{itemize}
\end{lemma}
\begin{proof}
  Clearly $A_G$ satisfies \AX{t1} and \AX{t3}. If $A_G(u,w)\ne\emptyset$
  and $A_G(w,v)\ne\emptyset$ there is a $u,w$-path $P_1$ and a $w,v$-path
  $P_2$.  Since $w\in V(P_1)\cap V(P_2)\ne\emptyset$, the first vertex $x$
  along $P_1$ also contained in $P_2$ is well-defined. Then the
  concatenation $Q$ of the $u,x$-subpath of $P_1$ and the $x,v$-subpath of
  $P_2$ is a $u,v$-path and thus $A_G(u,v)\ne\emptyset$, i.e., \AX{t0} is
  satisfied. In particular, therefore, $A_G$ is a directed transit
  function.

  Moreover, we have $x\in V(P_1)\subseteq A_G(u,w)$,
  $x\in V(P_2)\subseteq A_G(u,w)$, and $x\in V(Q)\subseteq A_G(u,v)$.  Thus
  $x\in A_G(u,w)\cap A_G(w,v)\cap A_G(u,v)\ne\emptyset$, and thus \AX{Mod}
  holds.

  If $x\in A_G(u,v)$ then in particular there exist a path from $u$ to $v$
  through $x$, which contains a sub-path from $u$ to $x$ and a sub-path
  from $x$ to $v$. Thus $A_G(u,x)\neq \emptyset$ and $A_G(x,v)\ne
  \emptyset$, i.e., \AX{tr2} is satisfied.
  
  Let $R(u,w)\cap R(w,v)=\{w\}$ and assume
  $R(u,w)\cup R(w,v)\not\subseteq R(u,v)$. Then there is an
  $x\in R(u,w)\cup R(w,v)$ such that $x\notin R(u,v)$. Assume $x\in R(u,w)$
  then there is a path from $u$ over $x$ to $w$. Since $x\notin A_G(u,v)$
  but there is at least one path from $w$ to $v$, the $uxw$-path and the
  $wv$-path must overlap. But we have $A_G(u,w)\cap A_G(w,v)=\{w\}$, a
  contradiction. An analogous argument can be made for the case that
  $x\in A_G(w,v)$. Thus $A_G$ satisfies \AX{b5}.

  If $R(u,w)\neq\emptyset$ and $R(w,v)\neq\emptyset$ then there are
  paths $P_1$ from $u$ to $w$ and $P_2$ from $w$ to $v$. Let $x$ be the
  first vertex along $P_1$ that is also contained in $P_2$. Denote by
  $P_1'$ the subpath of $P_1$ from $u$ to $x$ and by $P_2'$ the subpath of
  $P_2$ from $x$ to $v$. Then $P_1'\cup P_2'$ is a path from $u$ to $v$. 
\end{proof}
Lemma~\ref{lem:propdir} generalizes the results in \cite{ferdoos_Thesis},
where it is shown that the all-paths function $A_G$ of any connected
undirected graph satisfies \AX{b5} and \AX{Mod}.  Connectedness in this
setting implies that the precondition in \AX{Mod} is always satisfied.

\begin{figure}
  \begin{center}
    \begin{tikzpicture}
      \tikzset{vertex/.style = {shape=circle,draw,minimum size=1.0em}}
      \tikzset{edge/.style = {->,> = latex',dashed}}
      \node[vertex] (u) at  (0,0) {$u$};
      \node[vertex] (1) at  (3,0) {$x$};
      \node[vertex] (v) at  (6,0) {$v$};
      \node[vertex] (w) at  (3,-1.25) {$w$};
      
      \draw[edge] (u) to (1);
      \draw[edge] (1) to (v);      
      \draw[edge] (1) to [bend left=30] (w);
      \draw[edge] (w) to [bend left=30] (1);
    \end{tikzpicture}
  \end{center}
  \caption{Example of a directed graph that does not satisfy many of the
    axioms that are essential to prove the existence of paths for DAGs, or
    their posets. For example we do not have $R(u,w)\neq\emptyset$ and
    $R(w,v)\neq\emptyset$ implies that $w\in R(u,v)$, hence \AX{tr0} is not
    satisfied. Similarly, we have that $R(u,w)\cup R(w,v)\not\subseteq
    R(u,v)$ and \AX{tr1} does not hold. Since $R(u,w)\cap R(w,v)=\{w,x\}$
    also \AX{q} is not satisfied.}
\end{figure}
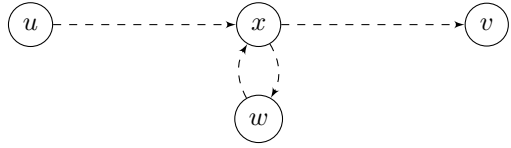

\begin{figure}[htb]
  \begin{tabular}{cc}
    \begin{subfigure}{0.45\textwidth}
      \begin{center}
	\begin{tikzpicture}
	  \tikzset{vertex/.style = {shape=circle,draw,minimum size=1.0em}}
	  \tikzset{edge/.style = {->,> = latex'}}
	  \node[vertex] (u) at  (0,0) {u};
	  \node[vertex] (x) at  (2,0) {x};
	  \node[vertex] (v) at  (4,0) {v};
	  \draw[edge] (u) to (x);	  
	  \draw[edge] (x) to[bend left=15] (v);
	  \draw[edge] (v) to[bend left=15] (x);
	  \draw[edge] (u) to[bend left=45] (v);
	\end{tikzpicture}
      \end{center}
      \subcaption{Counterexample \AX{b1$_1$}.}
    \end{subfigure}
    &
    \begin{subfigure}{0.45\textwidth}
      \begin{center}
	\begin{tikzpicture}
	  \tikzset{vertex/.style = {shape=circle,draw,minimum size=1.0em}}
	  \tikzset{edge/.style = {->,> = latex'}}
	  \node[vertex] (u) at  (0,0) {u};
	  \node[vertex] (x) at  (2,0) {x};
	  \node[vertex] (v) at  (4,0) {v};
	  \draw[edge] (u) to[bend left=15] (x);
	  \draw[edge] (x) to[bend left=15] (u);
	  \draw[edge] (x) to (v);
	  \draw[edge] (u) to[bend left=45] (v);
	\end{tikzpicture}
      \end{center}
      \subcaption{Counterexample \AX{b1$_2$}.}
    \end{subfigure}
    \\
    \begin{subfigure}{0.45\textwidth}
      \begin{center}
	\begin{tikzpicture}
	  \tikzset{vertex/.style = {shape=circle,draw,minimum size=1.0em}}
	  \tikzset{edge/.style = {->,> = latex'}}
	  \node[vertex] (u) at  (0,0) {u};
	  \node[vertex] (x) at  (2,0) {x};
	  \node[vertex] (v) at  (4,0) {v};
	  \node[vertex] (w) at  (3,-1) {w};
	  \draw[edge] (u) to (x);
	  \draw[edge] (x) to (v);
	  \draw[edge] (u) to[bend left=45] (v);
	  \draw[edge] (v) to (w);
	  \draw[edge] (w) to (x);
	\end{tikzpicture}
      \end{center}
      \subcaption{Counterexample \AX{b2$_1$}.}
    \end{subfigure}
    &
    \begin{subfigure}{0.45\textwidth}
      \begin{center}
	\begin{tikzpicture}
	  \tikzset{vertex/.style = {shape=circle,draw,minimum size=1.0em}}
	  \tikzset{edge/.style = {->,> = latex'}}
	  \node[vertex] (u) at  (0,0) {u};
	  \node[vertex] (x) at  (2,0) {x};
	  \node[vertex] (v) at  (4,0) {v};
	  \node[vertex] (w) at  (1,-1) {w};
	  \draw[edge] (u) to (x);
	  \draw[edge] (x) to (v);
	  \draw[edge] (u) to[bend left=45] (v);
	  \draw[edge] (x) to (w);
	  \draw[edge] (w) to (u);
	\end{tikzpicture}
      \end{center}
      \subcaption{Counterexample \AX{b2$_2$}.}
    \end{subfigure}\\
  \end{tabular}
  \caption{The all-path transit function $A_G$ of the graphs violates one
    of the betweenness axioms. (a) \AX{b1$_1$} is violated since $x\in
    R(u,v)=\{u,x,v\}$ but $v\in R(u,x)=\{u,v,x\}$. (b) \AX{b1$_2$} is
    violated by $x\in R(u,v)=\{u,x,v\}$ and $u\in R(x,v)=\{u,x,v\}$.  (c)
    \AX{b2$_1$} does not hold since $x\in R(u,v)=\{u,x,v\}$ and $w\in
    R(u,x)=\{u,v,w,x\}$ but $w\notin R(u,v)$ and hence $R(u,x)\not\subseteq
    R(u,v)$. (d) Analogously, \AX{b2$_2$} is violated by $x\in R(u,v)$ and
    $w\in R(x,v)$ but $w\notin R(u,v)$.}
  \label{fig:not_b1b2}
\end{figure}
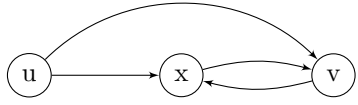
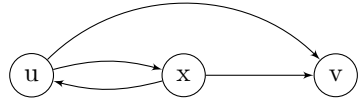
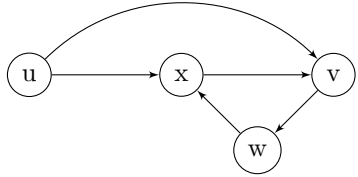
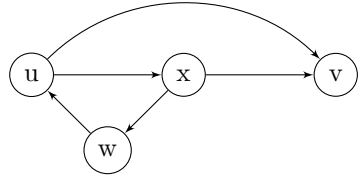

While the all-path transit function has many appealing properties in the
symmetric case, this is no longer the case for directed graphs.
Fig.~\ref{fig:not_b1b2} shows that $A_G$ in general satisfies neither
\AX{b1$_{1,2}$} nor \AX{b2}, and hence $A_G$ is not weakly geometric but
only satisfies the much weaker axiom \AX{b5}.  It is worth noting that
directed transit functions satisfying \AX{b4} and \AX{b5} also satisfy
\AX{b2}.  Moreover, \AX{b4} immediately implies (the contrapositive of)
  \AX{tr2}. Thus a directed transit function satisfying \AX{b4} and
  \AX{b5} is weakly geometric.

\begin{lemma}
  A directed transit function $R$ satisfying axiom \AX{b2} and
  \AX{Mod} also satisfies \AX{b5}.
\end{lemma}
\begin{proof}
  Suppose $R(u,x)\cap R(x,v)=\{x\}$.  By \AX{Mod}, we have $R(u,v)\cap
  R(u,x)\cap R(x,v)=R(u,v)\cap\{x\} \neq\emptyset$ and thus $x\in
  R(u,v)$. Now \AX{b2} immediately yield $R(u,x)\cup R(x,v)\subseteq
  R(u,v)$, and thus $R$ satisfies \AX{b5}.
\end{proof}
Note that a geometric directed transit function does not necessarily
satisfy \AX{b5}. As an example, consider $V=\{u,v,w\}$, $R(u,v)=\{u,v\}$,
$R(u,w)=\{u,w\}$, $R(w,v)=\{w,v\}$, $R(v,u)=R(w,u)=R(v,w)=\emptyset$,
and $R(x,x)=\{x\}$ for all $x\in V$. It is not difficult to check that $R$
is a geometric directed transit function. However, $R(u,w)\cap
R(w,v)=\{w\}$ but $R(u,v)=\{u,v\}$, violating \AX{b5}.

\subsection{Induced path transit function}

A path $P=(u=x_0,x_1,\dots,x_k=v)$ in a (directed) graph $G=(V,E)$ is
\emph{induced} if there is no edge $(i,j)\in E$ with $j\ne i-1$ for $1\le
j\le n$. That is, the subgraph $G[P]$ of $G$ induced by the vertex set $P$
comprises only the directed edges $(x_{i-1},x_i)$ for $1\le i\le k$.

A path $P$ is \emph{shortcut-free} if there is no edge $(i,j)\in E$ with
$j\ne i-1$ for $1\le j\le n$ and $i<j$. We caution the reader that what we
call here shortcut-free paths is sometimes called ``induced directed
paths'' \cite{Laetsch:10}.  For symmetric digraphs (and thus undirected
graphs) shortcut-free and induced paths coincide. In the directed case,
every induced path is shortcut-free, but the converse is not true.

The transit functions deriving from induced paths in undirected graphs have
been considered e.g.\ in \cite{Changat:99,Changat:04i,mcjmhm-10}.
Nebesk{\'y}'s proof that it is impossible to characterize the induced path
function of connected graphs using a set of first order axioms
\cite{Nebesky:02} suggests that it will be difficult to obtain general
results on this class of transit functions. Research in the undirected
cases therefore focuses on the characterization of induced path transit
functions for specific graph classes. In the directed case, it seems quite
natural to consider the shortcut-free paths.
\begin{definition}
  For a digraph $G$, let $x\in J_G(u,v)$ if $x$ lies along a shortcut-free
  path from $u$ to $v$. 
\end{definition}

\begin{lemma}
  The function $J_G$ is a directed transit function for every digraph $G$. 
\end{lemma}
\begin{proof}
  It follows immediately from the definition that $J_G$ satisfies \AX{t1}
  and \AX{t3}. Moreover, every paths $P$ in $G$ contains a shortcut-free
  subpath, obtained by removing, iteratively, from $P$ all vertices that
  are ``bridged'' by a shortcut. Thus $J_G(u,v)=\emptyset$ if and only if
  $A_G(u,v)=\emptyset$ for all $u,v\in V$. It follows immediately that
  $J_G$ satisfies \AX{t0} since $A_G$ satisfies \AX{t0}.
\end{proof}

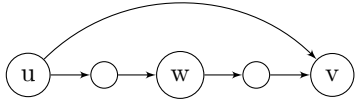
\begin{figure}
  \begin{center}
    \begin{tikzpicture}
      \tikzset{vertex/.style = {shape=circle,draw,minimum size=1.0em}}
      \tikzset{edge/.style = {->,> = latex'}}
      \node[vertex] (u) at  (0,0) {u};
      \node[vertex] (x) at  (1,0) {};
      \node[vertex] (w) at  (2,0) {w};
      \node[vertex] (y) at  (3,0) {};
      \node[vertex] (v) at (4,0)  {v};
      \draw[edge] (u) to (x);
      \draw[edge] (x) to (w);
      \draw[edge] (w) to (y);
      \draw[edge] (y) to (v);
      \draw[edge] (u) to [bend left=45] (v);
    \end{tikzpicture}
  \end{center}
  \caption{A graph $G$ for which the transit function $J_G$ does not
    satisfy \AX{b5}: we have $J_G(u,w)=\{u,...,w\}$,
    $J_G(u,w)=\{w,...,v\}$, but $J_G(u,v)=\{u,v\}$.}
  \label{fig:countercirc}
\end{figure}

The example in Figure~\ref{fig:countercirc} shows that $J_G$ in general
does not satisfy \AX{b5}.  The following two properties are satisfied by
The undirected induced path transit function satisfies the following two
properties, see e.g.\ \cite{mcjmhm-10}. Below, we show that this remains
true in the directed case:
\begin{lemma}
  The directed transit function $J_G$ of any digraph $G$ satisfies 
  \begin{itemize}
  \item[\AX{j2}] $R(u,x)=\{u,x\}$, $R(x,v)=\{x,v\}$, $R(u,v)\ne\{u,v\}$
    implies $x\in R(u,v)$.
  \item[\AX{j2'}] $x\in R(u,y)$ and $y\in R(x,v)$, $R(u,x)=\{u,x\}$,
    $R(x,y)=\{x,y\}$, $R(y,v)=\{y,v\}$, and $R(u,v)\ne\{u,v\}$
    implies $x\in R(u,v)$.
  \end{itemize}
\end{lemma}
\begin{proof}
  Let $G$ be a directed graph. If $u=v$ then $J_G(u,v)=\{u,v\}$, and thus
  the precondition of both \AX{j2} and \AX{j2'} is not met. Hence we may
  assume $u\ne v$. Setting $x=y$ in \AX{j2'} and removing the
  preconditions that become trivially true, i.e., $x\in R(u,x)$, $x\in
  R(x,v)$, and $R(x,x)=\{x,x\}$, shows that \AX{j2'} reduces to \AX{j2} for
  $x=y$. For $x=u$ or $x=v$ the precondition in \AX{j2} is never met, and
  hence the axiom holds trivially. If $u,v,x$ are pairwise distinct, then
  $(u,x)\in E(G)$ and $(x,v)\in E(G)$, $u\ne v$ and $(u,v)\notin E(G)$.
  Hence $(u,x,v)$ is a shortcut-free path of length $2$ in $G$. By
  definition, we therefore have $\{u,x,v\}\subseteq J_G(u,v)$ and thus in
  particular $x\in J_G(u,v)$, i.e., $J_G$ satisfies \AX{j2}.
  
  It thus remains to consider \AX{j2'} for $x\ne y$. Moreover, the
  precondition is void if $u=v$. For $x=u$ and $y=v$ or $x=v$ and $y=u$,
  the implication is trivially true.  Thus it remains to consider the case
  that $u,v,x,y$ are pairwise distinct. The precondition of \AX{j2'}
  implies that $(u,x)\in E(G)$, $(x,y)\in E(G)$, $(y,v)\in E(G)$ and thus
  $(u,x,y,v)$ is a directed path of length $3$.  Moreover, the condition
  $x\in R(u,y)$ requires that there is a shortcut-free path from $u$ to $y$
  that runs through $x$, and thus in particular $(u,y)\notin
  E(G)$. Analogously, $y\in R(x,v)$ implies $(x,v)\notin E(G)$. Finally,
  $R(u,v)\ne\{u,v\}$ implies $(u,v)\notin E(G)$, and thus $(u,x,y,v)$ is a
  shortcut-free directed path from $u$ to $v$ and thus
  $\{u,x,y,v\}\subseteq J_G(u,v)$.
\end{proof}

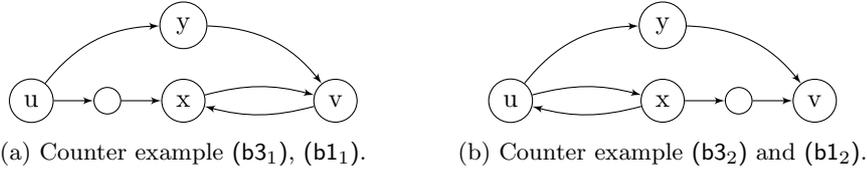
\begin{figure}[htb]
  \begin{tabular}{cc}
    \begin{subfigure}{0.45\textwidth}
      \begin{center}
        \begin{tikzpicture}
          \tikzset{vertex/.style = {shape=circle,draw,minimum size=1.0em}}
          \tikzset{edge/.style = {->,> = latex'}}
          \node[vertex] (u) at  (0,0) {u};
          \node[vertex] (z) at  (1,0) {};
          \node[vertex] (x) at  (2,0) {x};
          \node[vertex] (v) at  (4,0) {v};
          \node[vertex] (y) at (2,1)  {y};
          \draw[edge] (u) to (z);
          \draw[edge] (z) to (x);	  
          \draw[edge] (x) to[bend left=15] (v);
          \draw[edge] (v) to[bend left=15] (x);
          \draw[edge] (u) to[bend left=25] (y);
          \draw[edge] (y) to[bend left=25] (v);
        \end{tikzpicture}
      \end{center}
      \subcaption{Counter example \AX{b3$_1$}, \AX{b1$_1$}.}
    \end{subfigure}
    &
    \begin{subfigure}{0.45\textwidth}
      \begin{center}
        \begin{tikzpicture}
          \tikzset{vertex/.style = {shape=circle,draw,minimum size=1.0em}}
          \tikzset{edge/.style = {->,> = latex'}}
          \node[vertex] (u) at  (0,0) {u};
          
          \node[vertex] (x) at  (2,0) {x};
          \node[vertex] (z) at  (3,0) {};
          \node[vertex] (v) at  (4,0) {v};
          \node[vertex] (y) at (2,1)  {y};
          \draw[edge] (u) to[bend left=15] (x);
          \draw[edge] (x) to[bend left=15] (u);
          \draw[edge] (x) to (z);
          \draw[edge] (z) to (v);
          \draw[edge] (u) to[bend left=25] (y);
          \draw[edge] (y) to[bend left=25] (v);
        \end{tikzpicture}
      \end{center}
      \subcaption{Counter example \AX{b3$_2$} and \AX{b1$_2$}.}
    \end{subfigure}
  \end{tabular}
  \caption{Examples of graphs for which $J_G$ is not geometric/weakly geometric.}
  \label{J_G_not_geo_weakly}
\end{figure}

The induced path function $J_G(u,v)$ is not (weakly) geometric as shown by
the following examples Fig.~\ref{J_G_not_geo_weakly}.
  
\subsection{Directed interval functions}

A finite quasimetric space $(X,d)$ consists of a finite set $V$ and a
distance $d:V\times V\to \mathbb{R}_+\cup\{\infty\}$ satisfies, for all
$x,y,z\in V$, (i) $d(x,y)=0$ iff $x=y$ and (ii) the triangle inequality
$d(x,z)+d(z,y)\ge d(x,y)$.
\begin{definition}
  Let $(X,d)$ be a quasimetric space. Then the function $I_d:V\times V\to
  2^V$ given by
  \begin{equation}
    I_d(u,v) \coloneqq \begin{cases}
      \{ w\mid d(u,w)+d(w,v)=d(u,v) \} & \text{if } d(u,v)< \infty \\
      \emptyset                     & \text{if } d(u,v)= \infty \\
      \end{cases}
  \end{equation}
  is called the \emph{interval function} of $(X,d)$.
\end{definition}
This definition is the obvious analog of the undirected interval function
defined in \cite{Mulder:08a,Changat:01}. For strongly directed graphs, the
related set system $\{I_d(u,v)|u,v\in V(G)\}$ was used in \cite{Sur:08} to
define a convexity for digraph $G$.

The next result shows that the directed version shares the basic of
properties of the undirected interval function:
\begin{lemma}
  The interval function $I_d$ of a finite quasimetric space $(X,d)$
  is a geometric directed transit function.
\end{lemma}
\begin{proof}
  \AX{t1} follows from $d(u,v)=d(u,u)+d(u,v)=d(u,v)+d(v,v)$ since
  $d(u,u)=d(v,v)=0$ by (i), and \AX{t3} follow from
  $d(u,u)=0=d(u,w)+d(w,u)$ only if $w=u$.  If $I_d(u,w)\ne\emptyset$ and
  $I_d(w,v)\ne\emptyset$ then $d(u,w)<\infty$ and $d(w,v)<\infty$. The
  triangle inequality thus implies $d(u,v)\le d(u,w)+d(w,v)<\infty$ and
  thus $I_d(u,v)\ne\emptyset$, i.e., \AX{t0} is satisfied. Thus $I_d$ is
  a directed transit function.

  Now suppose $d(u,x)=\infty$ or $d(x,v)=\infty$. In the first case, we
  have $d(u,v)=\infty$ and thus $I_d(u,y)=\emptyset$, which in turn implies
  $x\notin I_d(u,v)$. Otherwise, $d(u,v)<\infty$, in which case
  $d(u,x)+d(x,v)\ne d(u,v)$ and thus $x\notin I_d(u,v)$. Hence $I_d$
  satisfies \AX{tr2}.

  Next consider the betweenness axioms.
  $x\in I_d(u,v)$ and $y\in I_d(u,x)$ implies $d(u,v)=d(u,x)+d(x,v)$ and
  $d(u,x)=d(u,y)+d(y,x)$, and thus $d(u,v)=d(u,y)+d(y,x)+d(x,v)$. Now
  suppose $d(y,x)+d(x,v)\ne d(y,v)$. Then the triangle inquality implies
  $d(y,x)<d(y,x)+d(x,v)$ and thus $d(u,v)>d(u,y)+d(y,v)$, violating the
  triangle inequality. Thus $d(y,x)+d(x,v)=d(y,v)$, i.e., $x\in I_d(y,v)$,
  i.e., \AX{(b3)$_1$} holds.  Moreover, we have $d(u,y)+d(y,v)=d(u,v)$,
  i.e., $y\in I_d(u,v)$ for all $y\in I_d(u,x)$, and thus
  $I_d(u,x)\subseteq I_d(u,v)$. An analogous argument yield \AX{b2$_2$} and
  \AX{b2$_2$}, and hence \AX{b2} is satisfied.
\end{proof}

It is not difficult to check that Lemma~\ref{lem:upathT-T} remains true for
geometric directed transit function. Recall that we have shown in
Theorem.~\ref{thm:connectedness} above that for a geometric directed
transit function, the underlying graph $G_R$ contains paths contained
entirely in $R(u,v)$ and running through any chosen $w\in R(u,v)$ provided
$R(u,v)\ne\emptyset$. In particular, therefore, there is $x,y\in R(u,v)$
such that $R(u,x) = \{u,x\}$ and $R(y,v) =\{y,v\}$.  To show that $x$ and
$y$ are unique if \AX{p} is assumed, it suffices to to appeal \AX{b4}
instead of \AX{q} in the proof Lemma~\ref{lem:upathT-T}.  It follows
immediately that the entire path from $u$ to $v$ is uniquely determined:
\begin{observation} 
  If $R$ is a geometric directed transit function that in addition
  satisfies \AX{p}, then $R(u,v)\ne\emptyset$ implies that there is a
  unique path from $u$ to $v$ in $G_R$.
\end{observation}

It is natural to ask under which conditions a geometric directed transit
function $R$ derives from a graph $G$ such that $R=I_d$ and $G_R=G$, where
$I_d$ is the interval function of the quasi-metric defined by the shortest
paths in $G$. For the symmetric case, i.e., for undirected graphs and their
metrics, several characterizations for such ``graphic'' transit functions
have been given \cite{nebe-94,Nebesky:98}. In the directed case, however,
this problem is more difficult and will be tackled in the forthcoming work.

\section{Concluding Remarks}

We have introduced here a notion of transit functions for directed and not
necessarily connected structures. It is a proper generalization of the
well-studied undirected transit functions for connected structures that
also encompasses the extension this framework to disconnected graphs in
\cite{Changat:17}. Directed transit functions also capture the idea of
directed betweenness that is common in network science \cite{White:94}. As
a consequence, we obtain a simple characterization of directed transit
functions that describe posets. This is in contrast to the difficulties of
describing posets by means of connected, undirected transit functions
\cite{brev2009,mathews2008,vandeVel:93}.

We observed that the most widely used concepts and in particular the
well-established betweenness axioms naturally generalize to the directed
case. For example, there is natural notion of weakly geometric transit
functions that ensures that $R(u,v)\ne\emptyset$ implies that there is a
directed path from $u$ to $v$ in the underlying digraph $G_R$, generalizing
a basic fact from the theory of undirected transit functions. Similarly, it
is possible to define many of the walk- and path-based constructions of
transit function on graphs to the directed case. We briefly consider the
directed analogs of interval (geodetic) functions, as well as all-paths,
and induced path transit functions. For directed acyclic graphs, we have
obtained characterizations closely matching results for undirected transit
functions on cycle-free graphs and trees. On the other hand, there are also
some quite drastic differences. While interval functions of digraphs turn
out to be a special case of geometric directed transit functions in
complete analogy with the undirected case, it does not seem possible to
extend the graphic transit functions \cite{nebe-94,Nebesky:98} to the
directed case. For further work on this topic we refer to forthcoming work.

Nebesk{\'y} \cite{nebe-02} showed that no first order characterization of
the induced-path function exists for undirected graphs. Since the symmetric
case is obtained from the definition of $J_G$ by simply adding the symmetry
axiom \AX{t2s}, such a characterization also may not exists for the the
general case. The all-paths transit function of undirected graphs has a
well-known characterization in terms of first order axioms
\cite{Changat:01}.  These conditions do not appear to generalize to
directed graphs, however, since it does not satisfy \AX{b2} for directed
graphs.  It remains an open problem, therefore, whether a simple
characterization for the all-paths directed transit functions can be
obtained.

Even if a simple characterization for transit functions deriving from
graphs is not possible, it remains a rich and interesting topic to identify
transit functions of particular graphs classes. In the undirected case, for
instance, there is a rich literature on distance hereditary graphs, see
e.g.\ \cite{Howorka:77,Bandelt:86}. The concept is extended to directed
graph in \cite{Laetsch:10}: A directed graph $G=(V,E)$ is \emph{distance
hereditary} if for all induced sub-digraphs $G'$ and all $u,v\in V(G')$ we
have $d_{G'}(u,v)=\infty$ or $d_{G'}(u,v)=d_G(u,v)$.  As shown in Lemma~1
of \cite{Laetsch:10}, $G$ is a distance hereditary if and only if every
shortcut-free directed path from $u$ to $v$ is a shortest directed path. In
\cite{Changat:24} it it show that the interval function of undirected
distance hereditary graph has a simple characterization in terms of first
order axioms. Similar results have been obtained for various graph classes
for interval function, for example see \cite{chalopin2024} and for induced
path see\cite{mcjmhm-10,Changat:19,Changat:23}. It will be interesting to
explore this type of results also of directed transit functions and
digraphs. Related questions arising from the present concern e.g.\ the
characterization of directed graphs $G$ in which the all-paths function
$A_G$ or the induced path function $J_G$ is (geometric/weakly geometric).
In recent years also the betwenness relations defined other undirected path
systems, such a triangle paths or toll walks have been considered
\cite{lcp,weaktoll}.We expect that some of these notions of betweenness
will have interesting generalizations to directed paths.

Transit functions also feature interesting connections to other set systems
and relations. Among symmetric transit functions, the geometric transit
functions can be characterized in terms of their ``base relations'' $\leq_b$
defined as
\begin{equation}
  \label{eq:eqvandeval}
  u\leq_b v \iff R(b,u)\subseteq R(b,v) \qquad\text{for all }u,v,b\in V
\end{equation}
which must be a partial order satisfying $R(u,v)=\{x\mid u \leq_b x\leq_b
v\}$, see \cite[Prop.~5.2]{vandeVel:93} and \cite{Bandelt:96}. In the
general case, however, equ.(\ref{eq:eqvandeval}) yields $u\leq_b v$ and
$v\leq_b u$ whenever $R(b,u)=R(b,v)=\emptyset$, i.e., $\leq_b$ is not a
total order on $V$.  It will be interesting to explore whether this
construction can be adjusted to yield a similar characterization for
geometric directed transit functions.

Taken together, the directed transit functions introduced here appear to be
a natural and useful generalization of the well-studied transit functions.

\begin{small}
\section*{Declarations}

% \subsection{Availability of Data and Materials}

\subsection*{Competing interests}
The authors declare that they have no competing interests, or other
interests that might be perceived to influence the results and/or
discussion reported in this paper.

\subsection*{Funding}
This work was supported in part by the DST, Govt. of India (Grant
No. DST/INT/DAAD/P-03/ 2023 (G)), the DAAD, Germany (Grant No. 57683501),
and the CSIR-HRDG for the Senior Research Fellowship
(09/0102(12336)/2021-EMR-I).

\subsection*{Authors' contributions}
MC, PGNS, and PFS designed the study. All authors contributed the
mathematical results and the drafting of the manuscript. 
% \NEW{All authors
%  participated in revising the manuscript during the review process.}
\end{small}
\bibliography{DTF}

\end{document}